%% file: _main.tex
\documentclass[a4paper,10pt]{elsarticle}
\usepackage{jcomp}
\usepackage{framed,multirow}
\usepackage[utf8]{inputenc}
\usepackage{bm}
\usepackage{amsmath}
\usepackage{mathrsfs}
\usepackage{graphicx}
\usepackage{epsfig, setspace}
\usepackage{url}
\usepackage{graphicx, transparent}
\usepackage{color}
\usepackage{algorithm} 
\usepackage{algorithmicx}
\usepackage{etex}
\usepackage{dirtytalk}
\usepackage[mathscr]{euscript}
\usepackage[bookmarks=true]{hyperref}
\usepackage{amsmath}
\usepackage{xcolor, soul}
\usepackage{rotating} 
\usepackage{pdflscape} 
\usepackage{silence}
\usepackage{ulem,cancel}
\newtheorem{remark}{\bf Remark}[section]

\usepackage{graphicx}
\usepackage{float}
\usepackage{subfigure}
\usepackage{subfigmat}

\usepackage{tikz}
\usepackage{capt-of}
\usepackage{setspace,lipsum}
\usepackage{lineno}
\usepackage{gensymb}
\usepackage{amssymb}
\biboptions{sort&compress}

\usepackage{setspace,lipsum}
\usepackage{listings}
\definecolor{aquamarine}{rgb}{0.5, 1.0, 0.83}
\definecolor{OliveGreen}{rgb}{0,0.6,0}
\usepackage{booktabs}
\usepackage{cellspace}
\usepackage{rotating}
\usepackage{tikz}
\usepackage{mathrsfs}

\sethlcolor{aquamarine}

\definecolor{codegreen}{rgb}{0,0.6,0}
\definecolor{codegray}{rgb}{0.5,0.5,0.5}
\definecolor{codepurple}{rgb}{0.58,0,0.82}
\definecolor{backcolour}{rgb}{0.95,0.95,0.92}

\lstdefinestyle{mystyle}{
    backgroundcolor=\color{backcolour},   
    commentstyle=\color{codegreen},
    keywordstyle=\color{magenta},
    numberstyle=\tiny\color{codegray},
    stringstyle=\color{codepurple},
    basicstyle=\ttfamily\footnotesize,
    breakatwhitespace=false,         
    breaklines=true,                 
    captionpos=b,                    
    keepspaces=true,                 
    numbers=left,                    
    numbersep=5pt,                  
    showspaces=false,                
    showstringspaces=false,
    showtabs=false,                  
    tabsize=2
}

\usepackage{xcolor} 
\let\oldtextcolor\textcolor
\renewcommand{\textcolor}[2]{\oldtextcolor{black}{#2}}
\let\oldcolor\color
\renewcommand{\color}[1]{\oldcolor{black}}

\begin{document}

\begin{frontmatter}
\title{High Resolution Optimized High-Order Schemes for Discretization of Non-Linear Straight and Mixed Second Derivative Terms}

\author[]{Hemanth Chandravamsi\footnote{Corresponding author, hemanthgrylls@gmail.com}}
\author[]{Steven H.\ Frankel}
\address{Faculty of Mechanical Engineering, Technion - Israel Institute of Technology, Haifa, Israel}


\begin{abstract}
In this paper, we propose a new set of midpoint-based high-order discretization schemes for computing straight and mixed nonlinear second derivative terms that appear in the compressible Navier-Stokes equations. Firstly, we detail a set of conventional fourth and sixth-order baseline schemes that utilize central midpoint derivatives for the calculation of second derivatives terms. To enhance the spectral properties of the baseline schemes, an optimization procedure is proposed that adjusts the order and truncation error of the midpoint derivative approximation while still constraining the same overall stencil width and scheme order. A new filter penalty term is introduced into the midpoint derivative calculation to help achieve high wavenumber accuracy and high-frequency damping in the mixed derivative discretization. Fourier analysis performed on the both straight and mixed second derivative terms show high spectral efficiency and minimal numerical viscosity with no odd-even decoupling effect. Numerical validation of the resulting optimized schemes is performed through various benchmark test cases assessing their theoretical order of accuracy and solution resolution. The results highlight that the present optimized schemes efficiently utilize the inherent viscosity of the governing equations to achieve improved simulation stability - a feature attributed to their superior spectral resolution in the high wavenumber range. The method is also tested and applied to non-uniform structured meshes in curvilinear coordinates, employing a supersonic impinging jet test case.

\end{abstract}


\end{frontmatter}

\input{introduction}

\input{governing_equations}
\input{methodology}

\input{timeDisc}

\input{results}

\input{conclusions}

\input{acknowledgements}
\input{appendix}

\bibliographystyle{elsarticle-num}
\bibliography{curvilinear}

\end{document}

%% file: introduction.tex

\tableofcontents

\section{Introduction}\label{sec:intro}

One of the crucial factors that has enabled present-day scale-resolving fluid flow simulations is the development of high-resolution numerical schemes \cite{wang2013high, laizet2009high, delorme2021application}. In this regard, a great deal of research has been dedicated to improving the resolution capabilities of both inviscid \cite{lele1992compact, bogey2004family, liu2013new, shu1999high} and viscous flux discretization \cite{liu2009direct, huynh2009reconstruction, gaitonde1998high, shen2010large, hermeline2009finite, de2021high}. While inviscid fluxes contribute to the generation of broadband turbulent features and characteristic waves in the flow, viscous fluxes, on the other hand, add dissipation and dictate the range of spatial flow scales that should be present in the flow at a given Reynolds number. The need for the development of low dispersion and low dissipation error numerical schemes is necessary to appropriately model inviscid and viscous effects. \\


The two important properties typically considered while designing a numerical scheme are the order of accuracy and spectral resolution (amongst other include stability and consistency). Firstly, the order of accuracy is an important consideration, as it dictates the response of numerical error to the change in grid resolution. However, order of accuracy alone cannot ensure an accurate scheme. The spectral properties of the numerical scheme is arguably an even more important factor, as they determine how accurately the scheme can capture various Fourier modes within the solution on a specific grid. \textcolor{red}{In other words, spectral properties dictate how discretization error responds to changes in the solution wavenumber on a fixed grid.} This directly determines the grid count needed to resolve a wave of particular frequency. This has lead to the birth to the class of dispersion relation preserving (DRP) schemes designed for first derivatives \cite{tam1993dispersion}. In the context of designing diffusion schemes, the typical focus in the literature is generally being devoted on the order of accuracy \cite{shen2009high, shen2010large, wang2018accuracy,de2021high} and the analysis of the scheme in the spectral space is often not considered. As pointed out by Nishikawa \cite{nishikawa2011robust, nishikawa2017effects, nishikawa2010beyond, nishikawa2011two} in the context of unstructured grids, the $h$-elliptic property (high frequency damping) of the diffusion scheme can be an important factor for solution convergence and robustness. Moreover, considering the inclusion of both straight and mixed derivative terms in the viscous fluxes of the compressible Navier-Stokes equations, it becomes necessary to handle the discretization of these terms separately. To the best of our knowledge, there are no reported studies in the literature that discuss the discretization of mixed derivative terms (e.g., of the form $\frac{\partial}{\partial x} . \frac{\partial \phi}{\partial y}$) in detail that include spectral analysis. These considerations motivate the development of high-resolution numerical schemes for the discretization of both straight and mixed second derivative terms that find application in a broad range of partial differential equations, including but not limited to the compressible Navier-Stokes equations, nonlinear diffusion equation, and Burgers' equations.





\subsection{Review of viscous flux discretization schemes for compressible Navier-Stokes Equations}

Various approaches in the literature for discretizing compressible viscous terms can be divided into five main categories. To motivate the present work, we discuss each of these methods here briefly. In a simplified context, the straight diffusion terms within the Navier-Stokes equations can be expressed using the following one-dimensional diffusion term.

\begin{equation} \label{straight}
    \frac{\partial f}{\partial x}=\frac{\partial}{\partial x}\left( \nu \frac{\partial \phi}{\partial x}\right),
\end{equation}

\noindent where, $f$ represents the diffusion flux, $\nu$ represents the diffusivity coefficient, and $\phi$ represents the field variable, i.e., velocity or temperature that appear in the Navier-Stokes equations. The following are various approaches in the literature used to discretize this term:

\begin{enumerate}
    \item \textbf{Successive application of two first derivatives} (Visbal and Gaitonde \cite{Visbal2002}):\\
    This approach which is used widely \cite{fang2014investigation,visbal1999high} requires computing gradients at all nodal locations twice: once for the gradients of the field variable $\phi$ and again for the flux $f$. A second-order version of this approach can be expressed as follows:

    \begin{equation}
        \left( \frac{\partial f}{\partial x} \right)_j \approx \frac{f_{j+1}-f_{j-1}}{2 \Delta x} \approx \frac{\nu_{j+1} \left( \frac{\phi_{j+2}-\phi_j}{2\Delta x}\right) - \nu_{j-1}\left( \frac{\phi_{j}-\phi_{j-2}}{2\Delta x}\right)}{2 \Delta x},
    \end{equation}

    With a constant diffusivity value $\nu$, the expression reduces to: $\nu \left( \phi_{j+2} - \phi_j + \phi_{j-2} \right)/\Delta x^2$, which only contains even indices $j+2$, $j$, and $j-2$. The lack of odd indices in the expressions can lead to the \textit{odd-even decoupling effect}. However, the approach of Visbal and Gaitonde \cite{Visbal2002} employs filtering to provide damping in the high-wavenumber range which suppresses the grid-to-grid modes effectively. It can also be noted that this approach does not guarantee conservation of flux `$f$'.

    \item \textbf{Derivative splitting through chain rule} (Gaitonde and Visbal \cite{gaitonde1998high}, Sandham and Li \cite{sandham2002entropy}):\\
    Gaitonde \& Visbal \cite{gaitonde1998high} have suggested that the flux derivative can be analytically split using the chain rule as follows:
    \begin{equation}
       \frac{\partial f}{\partial x}=\frac{\partial}{\partial x}\left( \nu \frac{\partial \phi}{\partial x}\right) = \nu \frac{\partial^2 \phi}{\partial x^2} + \frac{\partial \phi}{\partial x} \frac{\partial \nu}{\partial x}.
    \end{equation}

    This splitting approach was later comprehensively detailed and demonstrated by Sandham \& Li in \cite{sandham2002entropy}. A second order discretization of these split terms can be done as follows:

    \begin{equation}
       \left( \frac{\partial f}{\partial x} \right)_j \approx \nu \left( \frac{\phi_{j+1}-2\phi_j-\phi_{j-1}}{\Delta x^2}\right) + \left( \frac{\phi_{j+1}-\phi_{j-1}}{\Delta x} \right) \left( \frac{\nu_{j+1}-\nu_{j-1}}{\Delta x} \right) 
    \end{equation}

    The key advantage of this approach is that, any tendency of odd-even decoupling in the flow is partly suppressed by the diffusion terms with the help of the natural viscosity of the fluid. However, as specified in \cite{gaitonde1998high}, it is important to note that in three dimensions, particularly when solving in curvilinear coordinates and dealing with mixed derivatives, the number of gradients required to estimate various diffusion terms increases. Consequently, this leads to higher memory and CPU usage, resulting in reduced overall computational efficiency compared to other methods. Additionally, this approach does not hold the conservation property.
    
    \item \textbf{Conservative flux derivative of midpoint fluxes-1}  (Gaitonde and Visbal \cite{gaitonde1998high,visbal1999high}, Vanna et al. \cite{de2021high}):\\
    As opposed to computing derivatives of field variable $\frac{\partial \phi}{\partial x}$ at nodal locations $i$, the monograph of Gaitonde and Visbal \cite{gaitonde1998high} suggests the computation of derivatives at midpoint locations $f_{i+\frac{1}{2}}$ and use them to compute the flux derivatives in the next step as follows: 
    \begin{equation}
        \left( \frac{\partial f}{\partial x} \right)_j \approx \frac{ \nu_{j+\frac{1}{2}} \left( \frac{\partial \phi}{\partial x} \right)_{j+\frac{1}{2}} - \nu_{j-\frac{1}{2}} \left( \frac{\partial \phi}{\partial x} \right)_{j-\frac{1}{2}}}{\Delta x} \approx \frac{\nu_{j+\frac{1}{2}} \left( \frac{\phi_{j+1}-\phi_j}{\Delta x}\right) - \nu_{j-\frac{1}{2}} \left( \frac{\phi_{j}-\phi_{j-1}}{\Delta x}\right)}{\Delta x} 
    \end{equation}

    Vanna et al. in \cite{de2021high} have extended this approach to multi-fluid variable viscosity flows. Although their approach guarantees conservation, the scheme can only be at most second order accurate for non-uniform distributions of $\nu$ as demonstrated in \cite{de2021high}. However, as suggested in \cite{gaitonde1998high} high-order flux differentiation can be used to potentially achieve high-order accuracy which is one of the primary motivations for the present work.

    \item \textbf{Conservative flux derivative of midpoint fluxes-2} (Shen and Zha \cite{shen2009high,shen2010large}):\\
    The approach of Shen \& Zha \cite{shen2010large} also involves computing fluxes at midpoint locations. However, their stencil structure to compute flux `$f$' changes from one midpoint location to the other in the sense that some stencils have symmetric distribution of points about midpoint, while other have an asymmetric distribution, making use of most of the information from the stencil (except for the end nodes). This allows them to maintain a lower overall stencil width for the entire scheme while still achieving high-order accuracy. Although spectral analysis of their schemes is not presented anywhere in the literature (to the best of our knowledge), we note in the present work that their approach is also free from odd-even decoupling and provides high-frequency damping. Due to the extensive nature of the overall discretization procedure for their scheme, we have omitted the scheme's detailed explanation here for brevity. For more information, we refer interested readers to the original works \cite{shen2009high} and \cite{shen2010large}.
    
    \item \textbf{$\alpha$-damping approach} (Nishikawa \cite{Nishikawa2010} and Chamarthi et al. \cite{chamarthi2022importance}):\\
    The approach proposed by Nishikawa \cite{Nishikawa2010} utilizes cell-center values ($\phi$) and their respective cell-center derivatives ($\frac{\partial \phi}{\partial x}$) to estimate $\left(\frac{\partial f}{\partial x}\right)_j$ in a conservative fashion. The discretization as presented in \cite{Nishikawa2010} in the context of 1-D diffusion equation is as follows (terms from Ref. \cite{Nishikawa2010} are rearranged here):\\
    \begin{equation} \label{nishikawaDamp}
        \begin{aligned}
            \left( \frac{\partial f}{\partial x} \right)_j &= \frac{f_{j+\frac{1}{2}} - f_{j-\frac{1}{2}}}{\Delta x} + \mathcal{O}(\Delta x^2) \\
            &= \nu \left[ \frac{\alpha}{2} \left( \frac{\phi_{j+1}-2 \phi_j+\phi_{j-1}}{\Delta x^2} \right) + \left( 1-\frac{\alpha}{2} \right) \left( \frac{\phi_{j+2}-2 \phi_j+\phi_{j-2}}{4 \Delta x^2} \right) \right] + \mathcal{O}(\Delta x^2),
        \end{aligned}
    \end{equation}

    where,

    \begin{equation*}
        \begin{aligned}
            & f_{j+\frac{1}{2}} = \frac{\nu}{2} \left[ \left( \frac{\partial \phi}{\partial x}\right)_j + \left( \frac{\partial \phi}{\partial x}\right)_{j+1} \right] + \frac{\nu \alpha}{\Delta x} \left( \phi_R - \phi_L \right), \quad \left( \frac{\partial \phi}{\partial x}\right)_j = \frac{\phi_{j+1}-\phi_{j-1}}{2 \Delta x},\\
            &\left( \frac{\partial \phi}{\partial x}\right)_{j+1} = \frac{\phi_{j+2}-\phi_{j}}{2 \Delta x}, \quad \phi_L = \phi_j + \frac{\Delta x}{2} \left( \frac{\partial \phi}{\partial x} \right)_j, \quad \phi_R = \phi_{j+1} - \frac{\Delta x}{2} \left( \frac{\partial \phi}{\partial x} \right)_{j+1}.
        \end{aligned}
    \end{equation*}

    The diffusivity coefficient $\nu$ is treated constant here as presented in the original work. It is worth noting that the final expression in Eqn. \ref{nishikawaDamp} encompasses all values of $\phi$ within the index range of $[-2,2]$, effectively avoiding odd-even decoupling for all $\alpha \ne 0$. The parameter $\alpha$ primarily governs the damping characteristics of the scheme and yields fourth-order accuracy when set to $\alpha=8/3$. Subsequently, Chamarthi et al. \cite{chamarthi2022importance,chamarthi2023role} extended this approach to sixth-order and also highlighted the significance of high-frequency damping on the flowfield and the overall solution's accuracy. The respective schemes however degenerate to second order accuracy under variable viscosity (non-linear) conditions.
\end{enumerate}


While each of the above discussed methods has their own advantages, they also come with their respective drawbacks. For instance, although method-3 (Conservative flux derivative of midpoint fluxes-1) is conservative and provides high-frequency damping, its order of accuracy declines to second order, and it entails a large stencil size. Similarly, method-4, despite its compact stencil and $h$-ellipticity for straight derivative terms, it lacks resolution in the high-wavenumber region for the mixed derivative discretization \cite{shen2010large}. In addition, the stencil \cite{shen2010large} used in estimating each midpoint derivative does not incorporate information from all the points of the global second derivative stencil, resulting in a slightly reduced spectral resolution compared to what could be achieved. The objective of the current work is to bridge these gaps by developing a scheme that is (a) high-order accurate for linear and nonlinear cases, (b) conservative, and (c) high-frequency damping for both straight and mixed second derivative terms. Additionally, we aim to reduce the stencil size and employ an optimization process to enhance the spectral properties of the schemes, which will impact solution accuracy. The optimization of spectral properties and respective truncation errors have been previously explored for several inviscid flux discretization schemes \cite{tam1993dispersion,cheong2001grid,fu2017targeted,lin2018optimization,martin2006bandwidth,jin2018optimized,li2013optimized,ashcroft2003optimized}. However, to the best of our knowledge, there has been no study in the literature that focuses on the optimization of the spectral properties of viscous flux discretization scheme. Following are the various features and advantages of the class of schemes being proposed in the current work:

\begin{enumerate}
    \item \textbf{High-order:} Provides high-order accuracy for both linear and non-linear cases i.e., constant and variable viscosity flows.
    \item \textbf{Optimizability:} The schemes to be discussed in the present work offer flexibility in terms of tuning the spectral properties, while maintaining the stencil size and order of accuracy unchanged. This provides room for optimization. Examples of other such applications can be found in Refs. \cite{lamballais2011straightforward, lamballais2021viscous, fu2021very, dairay2017numerical}.
    \item \textbf{Spectral like resolution:} The optimized versions of the present schemes feature superior spectral accuracy in the high wavenumber range, resulting in accurate viscous dissipation. This leads to good solution accuracy and can potentially reduce the grid requirements for applications such as Direct Numerical Simulations (DNS).
    \item \textbf{Odd-even decoupling free:} The schemes to be presented are odd-even decoupling free and can suppress unwanted oscillations that may arise from the dispersion error of the inviscid schemes, thereby enhancing the simulation stability.
\end{enumerate}

\textcolor{magenta}{Of particular importance to the present discussion on viscous flux discretization are shear-dominated vortical flows. In such flows, high velocity gradients mean that dissipation from the second derivative terms in the governing equations potentially plays an important role. It is well established in the literature that instability growth in these flows is sensitive to various parameters, such as grid resolution/type, discretization approach, artificial/numerical dissipation, and the physical model employed (governing equations). More specifically, the impact of inviscid and viscous flux discretization on shear-dominated vortical flows, particularly at low Mach number settings, has been studied both theoretically and numerically \cite{tsoutsanis2015comparison, drikakis2001spurious}. Mosedale and Drikakis \cite{mosedale2007assessment} demonstrated the superiority of high-order inviscid discretization in resolving shear instability over low-order methods, which can lead to less accurate results. For flows containing multi-component discontinuities, Thornber and Drikakis \cite{thornber2008numerical} presented a modified Roe scheme to control numerical dissipation and accurately model instability growth. In their subsequent work, Thornber et al. \cite{thornber2008entropy} derived analytical formulae to model the entropy dissipation rate, demonstrating how excess dissipation can lead to inaccurate results in the context of Godunov-type nonlinear numerical discretizations, contrary to the conclusions of Minion and Brown \cite{minion1997performance}. We include two shear dominated vortical flow test cases in our numerical tests and discussion (section \ref{sec:results}) to study the effect of the present viscous flux discretization schemes: the doubly periodic shear layer and the Kelvin-Helmholtz instability.} \\





The remainder of the paper is structured as follows: the next section presents the governing equations, followed by the details of the proposed second derivative discretization schemes in Section \ref{sec:Viscmethod}. Section \ref{sec:Viscmethod} also focuses on the optimization and spectral error analysis of these schemes. The accuracy and the damping properties are investigated as well. Section \ref{sec:time-int} presents time integration and stability criterion. Section \ref{sec:results} presents numerical examples to demonstrate order of accuracy and solution resolution. Conclusions are presented towards the end.



%% file: governing_equations.tex
\section{Governing equations}   \label{sec:gov-eqns}

The non-dimensionalized 3-D compressible Navier-Stokes equations are being solved in the current work in Cartesian coordinates. The equations can be expressed in their conservative vector form as follows:

\begin{equation} \label{eqn:cns}
    \frac{\partial \mathbf{U}}{\partial t} + \frac{\partial \mathbf{F}^c}{\partial x} + \frac{\partial \mathbf{G}^c}{\partial y} + \frac{\partial \mathbf{H}^c}{\partial z} = \frac{1}{\mathrm{Re}} \left( \frac{\partial \mathbf{F}^v}{\partial x} + \frac{\partial \mathbf{G}^v}{\partial y} + \frac{\partial \mathbf{H}^v}{\partial z} \right),
\end{equation}

\noindent where $t$ represents time and $(x, y, z)$ denote the spatial coordinates. The conservative variable vector is $\mathbf{U} = \left[ \rho, \rho u, \rho v, \rho w, \rho E \right]^\top$. The inviscid fluxes $\mathbf{F}^c$, $\mathbf{G}^c$, and $\mathbf{H}^c$ are:

\begin{subequations}\label{eqn:invFluxes}
    \begin{align}
        \mathbf{F}^c &= \left[ \rho u, \rho u^2 + p, \rho u v, \rho u w, u \left(E+ p\right) \right]^\top, \\
        \mathbf{G}^c &= \left[ \rho v, \rho v u, \rho v^2 + p, \rho v w, v \left(E+ p\right) \right]^\top, \\
        \mathbf{H}^c &= \left[ \rho w, \rho w u, \rho w v, \rho w^2 + p, w \left(E+ p\right) \right]^\top,
    \end{align}
\end{subequations}

\noindent where $\rho$ is the density, $p$ is the pressure, and the velocity vector is $\Bar{\upsilon}=[u,v,w]$. The internal energy per unit volume is defined as $E = \frac{p}{\gamma-1} + \frac{\rho}{2}\left(u^2 + v^2 + w^2 \right)$. The equation of state $p \gamma \mathrm{M}^2= \rho T$ provides closure to the governing equations. The fluid is assumed to be calorically perfect with a specific heat ratio of $\gamma=1.4$. The viscous flux vectors $\mathbf{F}^v$, $\mathbf{G}^v$, and $\mathbf{H}^v$ are:

\begin{subequations}\label{eqn:viscFluxes}
    \begin{align}
        \mathbf{F}^v=\left[0, \tau_{x x}, \tau_{x y}, \tau_{x z}, u \tau_{x x}+v \tau_{x y}+w \tau_{x z}-q_{x}\right]^\top, \\
        \mathbf{G}^v=\left[0, \tau_{x y}, \tau_{y y}, \tau_{y z}, u \tau_{y x}+v \tau_{y y}+w \tau_{y z}-q_{y}\right]^\top, \\
        \mathbf{H}^v=\left[0, \tau_{x z}, \tau_{y z}, \tau_{z z}, u \tau_{z x}+v \tau_{z y}+w \tau_{z z}-q_{z}\right]^\top.
    \end{align}
\end{subequations}

\noindent The viscous normal and shear stresses in accordance to Newton's law of viscosity can be expressed as:

\begin{subequations}\label{visc}
    \begin{gather}
        \tau_{xx} = 2 \mu \frac{\partial u}{\partial x} + \lambda \left( \mathbf{\nabla} \cdot \Bar{\upsilon} \right),
        \tag{\theequation a-\theequation c}
        \quad
        \tau_{yy} = 2 \mu \frac{\partial v}{\partial y} + \lambda \left( \mathbf{\nabla} \cdot \Bar{\upsilon} \right), \quad
        \tau_{zz} = 2 \mu \frac{\partial w}{\partial z} + \lambda \left( \mathbf{\nabla} \cdot \Bar{\upsilon} \right),\\
        \tau_{xy} = \tau_{yx} = \mu \left(\frac{\partial u}{\partial y} + \frac{\partial v}{\partial x} \right),
        \quad
        \tau_{yz} = \tau_{zy} = \mu \left(\frac{\partial v}{\partial z} + \frac{\partial w}{\partial y} \right),
        \quad
        \tau_{xz} = \tau_{zx} = \mu \left(\frac{\partial u}{\partial z} + \frac{\partial w}{\partial x} \right),
        \tag{\theequation d-\theequation f}
    \end{gather}
\end{subequations}

\noindent where $\mu$ is the non-dimensional temperature dependent dynamic viscosity computed using the Sutherland's law \cite{Sutherland1893}. As a result of non-dimensionalization, the parameters Reynolds number ($\mathrm{Re}$) and Mach number ($\mathrm{M}$) are obtained. These quantities are defined based on reference values as $\mathrm{Re} = (\rho_{\text{ref}} L_{\text{ref}} u_{\text{ref}})/\mu_{\text{ref}}$ and $\mathrm{M}=u_{\text{ref}}/\sqrt{\gamma R T_{\text{ref}}}$. $\lambda = -\frac{2}{3} \mu$ is considered with Stokes' hypothesis \cite{kundu2015fluid}. The heat diffusion is modeled through Fourier's law of heat conduction as follows:

\begin{subequations}
    \begin{gather}
        q_{x} = -\kappa \frac{\partial T}{\partial x},
        \quad
        q_{y} = -\kappa \frac{\partial T}{\partial y},
        \quad
        q_{z} = -\kappa \frac{\partial T}{\partial z},
        \tag{\theequation a-\theequation c}
    \end{gather}
\end{subequations}

\noindent where $\kappa = \mu / (\mathrm{M}^2 (\gamma-1)\mathrm{Pr})$ is the scaled thermal conductivity after non-dimensionalization and $T$ denotes temperature. The inviscid discretization scheme and the proposed viscous flux discretization schemes to solve these equations are discussed next. 

%% file: methodology.tex
\section{Methodology} \label{sec:Viscmethod}
A new set of fourth- and sixth-order finite difference schemes for viscous flux discretization is presented in this section. For simplicity, the discretization method is discussed using the term $\frac{\partial \tau_{xx}}{\partial x}$ as an example from the $x$-momentum equation (Eqn. \ref{eqn:viscFluxes}). Extension to other viscous terms and directions will involve the same procedure. Consider the following expression:

\begin{equation} \label{tauxxTerm}
    \frac{\partial \tau_{xx}}{\partial x} = \frac{\partial }{\partial x} \left[\mu \left( \frac{4}{3} \frac{\partial u}{\partial x} - \frac{2}{3} \frac{\partial v}{\partial y} - \frac{2}{3} \frac{\partial w}{\partial z} \right) \right]
\end{equation}

The conservative flux derivative of $\tau_{xx}$ can be evaluated numerically through the following general expression:
\begin{equation} \label{eqn:2ndDer}
    \left. \frac{\partial \tau_{xx}}{\partial x}\right|_j = \frac{a^*}{\Delta x} \left[ \left(\tau_{xx}\right)_{j+\frac{1}{2}} - \left(\tau_{xx}\right)_{j-\frac{1}{2}}  \right] + \frac{b^*}{3 \Delta x} \left[ \left(\tau_{xx}\right)_{j+\frac{3}{2}} - \left(\tau_{xx}\right)_{j-\frac{3}{2}}  \right] + \frac{c^*}{5 \Delta x} \left[ \left(\tau_{xx}\right)_{j+\frac{5}{2}} - \left(\tau_{xx}\right)_{j-\frac{5}{2}}  \right] + TE.
\end{equation}



By appropriately selecting the coefficients $a^*$, $b^*$, and $c^*$, and employing the discretization method to evaluate the midpoint stresses such as $(\tau_{xx})_{j+\frac{1}{2}}$, a desired truncation error (TE) can be obtained. Since the shear stress terms involve velocity derivatives in all three spatial directions, a genuinely high-order accurate scheme can only be derived if and only if the following truncation errors are maintained while discretizing $\frac{\partial \tau_{xx}}{\partial x}$.

\begin{equation} \label{TEs}
    \begin{aligned}
        TE &= \mathcal{O}[\Delta x^4, \Delta y^4, \Delta z^4] \quad \longrightarrow \text{for 4th order schemes,} \\
        TE &= \mathcal{O}[\Delta x^6, \Delta y^6, \Delta z^6] \quad \longrightarrow \text{for 6th order schemes.}
    \end{aligned}
\end{equation}

Given the shear stresses are evaluated using a high order discretization, performing a Taylor series expansion on the flux derivative $\frac{\partial \tau_{xx}}{\partial x}$ (as previously presented in \cite{lele1992compact, gaitonde1998high}), leads to the following coefficient values: $[a^*, b^*, c^*] = [9/8, -1/8, 0]$ for all fourth-order schemes, and $[a^*, b^*, c^*] = [75/64, -25/128, 3/128]$ for all sixth-order schemes. With $[a^*, b^*, c^*]$  being fixed, the properties of the viscous flux discretization scheme will depend solely on the discretization method used for evaluating the midpoint stresses. \\



It can be noted that in Eqn. \ref{tauxxTerm}, the presence of straight and mixed derivative terms necessitates a different discretization approach for each of the terms. For example, while differentiating $\frac{\partial \tau_{xx}}{\partial x}$, the terms $\mu \frac{4}{3} \frac{\partial u}{\partial x}$ and $-\mu \frac{2}{3} \frac{\partial v}{\partial y}$ within $\tau_{xx}$ give rise to straight and mixed second derivative terms, respectively. For the sake of clarity, we present the discretization of straight and mixed derivative terms separately in two different sections: Sec. \ref{sec:straight} and Sec. \ref{sec:mixedDers}. Furthermore, to assess the properties of the discretization schemes in spectral space, we initially introduce the discretization of the straight derivative term, assuming a constant dynamic viscosity $\mu$. Later, in Section \ref{sec:viscCoeff}, we outline the discretization procedure for cases involving variable viscosity. 



\subsection{Discretization of the straight second derivative terms} \label{sec:straight}


In this section, we present two sets of fourth and sixth-order schemes: one corresponding to the baseline and the other to the optimized approach (new). The initial baseline approach employs the conventional central difference schemes for approximating midpoint derivatives, following the approaches of Zingg et al. \cite{zingg2000comparison} and De Rango \cite{de1999aerodynamic}. However, to approximate the midpoint derivatives in optimized schemes, we change the stencil design in order to efficiently use the available information. It will be demonstrated that the optimized second derivative schemes, even with the use of a smaller stencil, can exhibit better spectral properties compared to the baseline schemes.

\subsubsection{Baseline 4th and 6th order second derivative schemes} \label{sec:baseStraight}
\noindent Consider Eqn. \ref{straight} from the introduction section but with constant diffusivity coefficient $\nu=1$.

\begin{equation}
    \frac{\partial f}{\partial x} = \frac{\partial }{\partial x}\left( \frac{\partial \phi}{\partial x} \right)
\end{equation} \\

Using the gradients at midpoints, the second derivative of the field variable $\phi$ (velocity or temperature) can be approximated as follows:

\begin{equation} \label{gen2ndDer}
    \left. \frac{\partial^2 \phi}{\partial x^2} \right|_j \approx \frac{a^*}{\Delta x} \left( \phi^{\prime}_{j+\frac{1}{2}} - \phi^{\prime}_{j-\frac{1}{2}}  \right) + \frac{b^*}{3 \Delta x} \left( \phi^{\prime}_{j+\frac{3}{2}} - \phi^{\prime}_{j-\frac{3}{2}}  \right) + \frac{c^*}{5 \Delta x} \left( \phi^{\prime}_{j+\frac{5}{2}} - \phi^{\prime}_{j-\frac{5}{2}}  \right)
\end{equation} \\

As specified in Eqn. \ref{eqn:2ndDer}, the coefficients $[a^*, b^*, c^*]$ assume the values $[9/8, -1/8, 0]$ for fourth-order accuracy, and $[75/64, -25/128, 3/128]$ for sixth-order accuracy. Similarly, the gradient of $\phi$ at the midpoint locations can be computed as follows \cite{zingg2000comparison,de1999aerodynamic}:

\begin{equation} \label{midDersBase}
    \left. \frac{\partial \phi}{\partial x} \right|_{j+\frac{1}{2}} \approx \frac{a^*_1}{\Delta x} \left( \phi_{j+1} - \phi_{j}  \right) + \frac{b^*_1}{3 \Delta x} \left( \phi_{j+2} - \phi_{j-1}  \right) + \frac{c^*_1}{5 \Delta x} \left( \phi_{j+3} - \phi_{j-2}  \right)
\end{equation} \\

Achieving the desired order of accuracy with coefficients $[a^*, b^*, c^*]$ is contingent upon computing the midpoint derivatives (as in Equation \ref{gen2ndDer}) with the same or higher order of accuracy. Therefore, mathematically, it follows that for the shortest possible overall stencil, $a^*=a^*_1$, $b^*=b^*_1$, and $c^*=c^*_1$. This leads to an overall nine-point stencil $[j-4, j-3, \ldots, j+4]$ for the fourth-order scheme and an eleven-point stencil $[j-5, j-4, \ldots, j+5]$ for the sixth-order scheme. With respect to the explicit nature of the midpoint derivatives used in the scheme, the fourth and sixth-order schemes are denoted as ME4-Base and ME6-Base, expanding to Midpoint-based Explicit 4th/6th order schemes base scheme.

\begin{remark}
    Examining Eqn. \ref{midDersBase}, it is evident that the sixth-order midpoint derivative $\phi_{j+\frac{1}{2}}$ only incorporates points within the range $j-2:j+3$. Unfortunately, it neglects the information corresponding to points $\phi_{j-5}, \phi_{j-4}, \ldots, \phi_{j-1}$, which are still part of the overall second derivative stencil. This under utilization makes this approach less efficient in terms of stencil compactness and information utilization. This motivates the use of asymmetric/partially one-sided stencils to compute the midpoint derivatives \cite{shen2009high}. This approach helps optimize the spectral properties of the scheme and reduce the stencil length by efficiently incorporating information within the prescribed stencil.
\end{remark}

Fourier error analysis is employed to validate the order of accuracy and assess the spectral resolution of the schemes in the wavenumber domain. In this analysis, we consider a series of Fourier modes (or pure harmonic functions) represented by $\mathcal{F}(x) = e^{ikx}$, where the wavenumber $k$ ranges from the smallest wavenumber $k=0$ to the grid cutoff wavenumber $k=\pi/2$. By using the discrete values of $e^{ikx}$ on a uniform grid, a numerical approximation of the second derivative is computed. Subsequently, the second derivative approximation of each Fourier mode with the exact analytical second derivative of the input Fourier mode $\mathcal{F}(x)$ is computed.\\

\begin{equation}
    \mathcal{F}(x)=e^{ikx}, \quad k=\frac{2 \pi}{L} n, \quad n=0,1,2, \ldots, N / 2
\end{equation}

\begin{equation}
    \mathcal{F}^{\prime \prime}(x)= -k^2e^{ikx} = -k^2\mathcal{F}(x)
\end{equation}

\noindent The modified wavenumber $k^{*}_{xx}$ for the exact case is:
\begin{equation}
    k^{*}_{xx} = \frac{\mathcal{F}^{\prime \prime}(x)}{\mathcal{F}(x)} = -k^2
\end{equation}

By substituting $\mathcal{F}(x)$ in to the ME4-Base difference scheme (Eqn. \ref{gen2ndDer}), one can obtain the following modified wavenumber expression:
\begin{equation}
    \begin{aligned}
        \left(k^{*}_{xx}\right)^{\text{ME4-Base}} &= -\frac{365}{144} + \frac{87 \cos (k)}{32}-\frac{3 \cos (2 k)}{16} +\frac{\cos (3 k)}{288} \\
        &=-k^2 \left(1-\frac{3 k^4}{320}+\frac{k^6}{1792}+\frac{k^8}{230400}+ \cdots \right)
    \end{aligned}
\end{equation} \label{ME4base-modified}

\noindent Similarly for the ME6-Base scheme,
\begin{equation}
    \begin{aligned}
        \left(k^{*}_{xx}\right)^{\text{ME6-Base}} &= -\frac{2539103}{921600} + \frac{12505 \cos (k)}{4096}-\frac{335 \cos (2 k)}{1024}+\frac{2245 \cos (3 k)}{73728}-\frac{5 \cos (4 k)}{4096}+\frac{9 \cos (5 k)}{204800} \\
        &=-k^2 \left(1-\frac{5 k^6}{3584}+\frac{25 k^8}{147456}-\frac{23 k^{10}}{2162688}+ \cdots \right)
    \end{aligned}
\end{equation} \label{ME6base-modified}

\begin{figure}[h!]
    \centering
    \includegraphics[width=70mm]{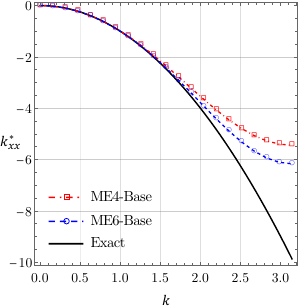}
    \caption{Modified wavenumber profiles of ME4-Base and ME6-Base schemes.}
    \label{baselineSchemes}
\end{figure}


\newpage
\subsubsection{Optimized 4th and 6th order second derivative schemes} \label{sec:optiStraight}

\begin{figure}[h!]
    \centering
    \includegraphics[width=0.85\textwidth]{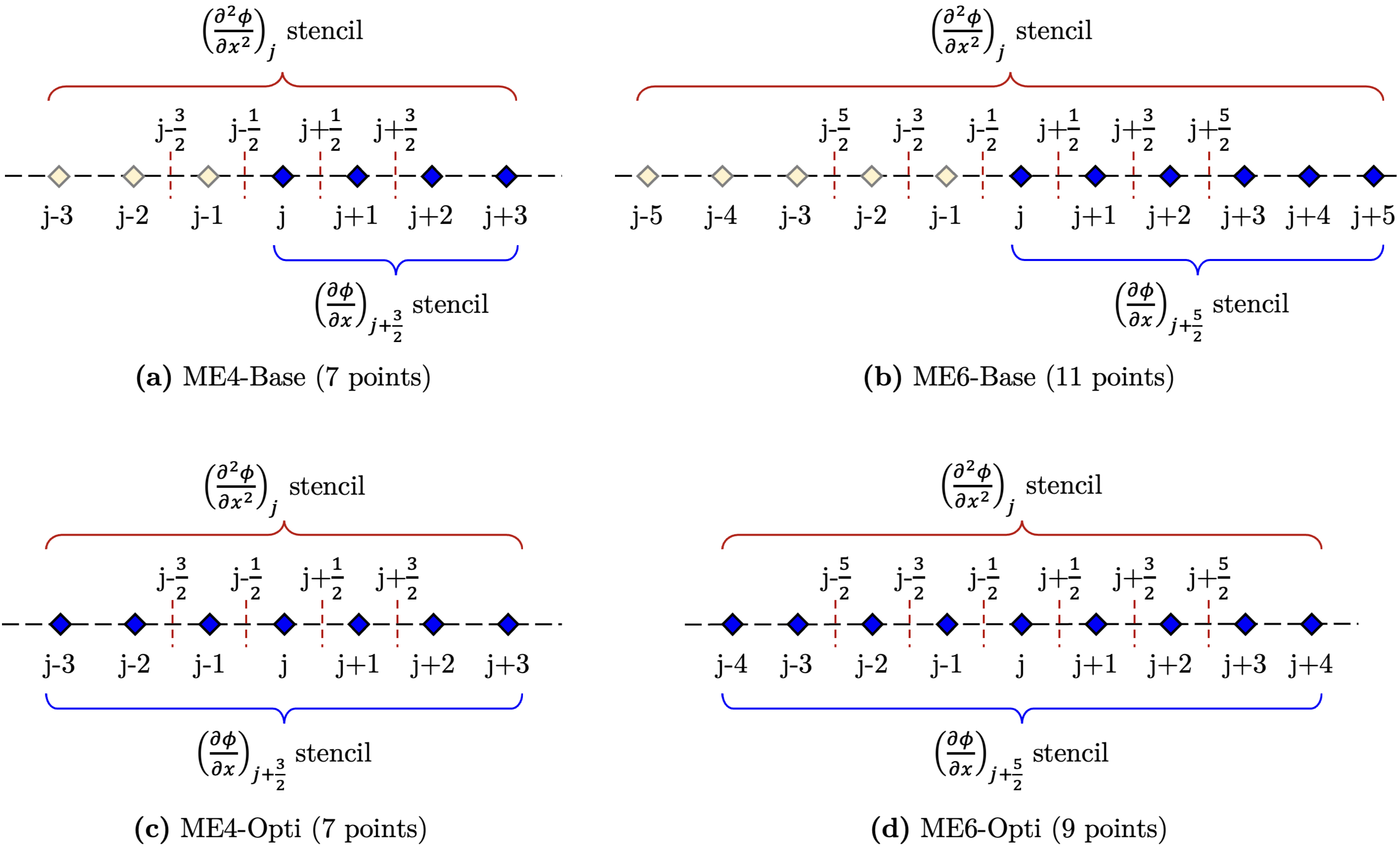}
    \caption{Stencil structures for the 4th and 6th order baseline/optimized second derivative schemes discussed in the present paper. The light yellow colored nodes represent the grid points that do not contribute towards the mid-point derivative calculation at $j+\frac{3}{2}$, and the blue points that do. The relatively larger stencil span of the midpoint derivatives corresponding to ME4-Opti and ME6-Opti schemes can be acknowledged.}
    \label{stencil}
\end{figure}

The optimized schemes presented in this section employ a seven point stencil for the fourth order scheme, and a reduced nine point stencil for the sixth order scheme in contrast to the ME6-Base scheme (which uses a eleven point stencil). The present optimized schemes will be shown to demonstrate improved spectral properties. This enhancement is primarily achieved by redesigning the stencil used for computing the midpoint derivatives. The present approach takes the advantage of \textit{all} the values inside the stencil to compute each midpoint derivative that is required for the calculation, which is not the case in the base schemes and other schemes from the literature \cite{shen2009high, shen2010large}. For instance, in the optimized fourth order scheme, the midpoint derivatives of $\frac{\partial \phi}{\partial x}$ at locations $j-\frac{3}{2}$, $j-\frac{1}{2}$, $j+\frac{1}{2}$, and $j+\frac{3}{2}$ all utilize a share stencil spanning from $j-3$ to $j+3$ (seven points), which is also the stencil for the final second derivative scheme. This is pictorially illustrated in Fig. \ref{stencil}(a,c). The same strategy has been employed for the optimized sixth-order scheme as well, which can be seen in Fig. \ref{stencil}(d). However, it utilizes a smaller stencil in comparison to the baseline scheme ME6-Base. This deliberate choice aims to minimize the number of floating-point operations, resulting in improved efficiency, as well as reduce the number of ghost cells at the boundaries. Following the nomenclature used for the baseline schemes, the optimized 4th and 6th order schemes presented in this section are referred to as ME4-opti and ME6-Opti respectively.

The general expressions for computing the midpoint derivatives used in constructing the fourth order ME4-Opti scheme, are given by the following Eqn. \ref{eqn:4thOrderOpti}.

\begin{equation}\label{eqn:4thOrderOpti}
    \begin{aligned}
    \left. \frac{\partial \phi}{\partial x}\right|_{j \pm \frac{1}{2}} & = \sum_{p=-3}^{3} a_{ \pm p} \frac{\pm\phi_{j+p}}{\Delta x}  \pm  \psi_{\frac{1}{2}} \frac{\partial^5\phi}{\partial x^5} \Delta x^5 + \mathcal{O}(\Delta x^6), \\
    \left. \frac{\partial \phi}{\partial x}\right|_{j \pm \frac{3}{2}} &= \sum_{p=-3}^{3} b_{ \pm p} \frac{\pm\phi_{j+p}}{\Delta x}  \pm  \psi_{\frac{3}{2}} \frac{\partial^5\phi}{\partial x^5} \Delta x^5 + \mathcal{O}(\Delta x^6).
    \end{aligned}
    \left\}
    \begin{aligned}
    \\
    \\
    \\
    \\
    \\
    \end{aligned}
    \right.
    \text{5th order}
\end{equation}

Where $\psi_{\frac{1}{2}}$ and $\psi_{\frac{3}{2}}$ are the coefficients of leading truncation error term in the Taylor series; an optimization procedure will be employed to find their values. \\

\begin{remark}
    In Eqn. \ref{eqn:4thOrderOpti}, although the utilization of a seven-point stencil enables a maximum accuracy of seventh order, the order of accuracy was intentionally limited to a lower value that is five. This restriction allows for control over the values of leading truncation error term coefficients $\psi_{\frac{1}{2}}$ and $\psi_{\frac{3}{2}}$, which facilitates the adjustment of final second derivative scheme's spectral properties; thus scope for improved spectral bandwidth resolution.
\end{remark}

Firstly, the coefficients $a_p$ and $b_p$ are determined as functions of $\psi_{\frac{1}{2}}$ and $\psi_{\frac{3}{2}}$, respectively. Expanding the expression $\left. \frac{\partial \phi}{\partial x}\right|_{j+\frac{1}{2}}$ in Eqn. \ref{eqn:4thOrderOpti}a around $j+\frac{1}{2}$ using Taylor series, and comparing the terms on left and right hand side, we obtain a family of coefficients that yield fifth-order accuracy for any given value of $\psi_{\frac{1}{2}}$. These coefficients can be expressed as follows:

\begin{equation} \label{apcoeffs}
\begin{aligned}
    a_{-3} = -\psi_{\frac{1}{2}}, \text{ } &a_{-2} = \frac{3}{640} (-1 + 1280 \psi_{\frac{1}{2}}), \text{ } a_{-1} = \frac{-5}{384} (-5 + 1152 \psi_{\frac{1}{2}}), \text{ } a_{0} = \frac{5}{64} (-15 + 256 \psi_{\frac{1}{2}}),\\
    a_{1} &= \frac{15}{64} (-5 + 64\psi_{\frac{1}{2}}), \text{ } a_{2} = \frac{1}{384} (-25 + 2304 \psi_{\frac{1}{2}}), \text{ } a_{-3} = \frac{1}{640} (3 - 640 \psi_{\frac{1}{2}})
\end{aligned}
\end{equation}

Similarly, for arbitrary $\psi_{\frac{3}{2}}$, a set of fifth-order coefficients for $b_p$ can be derived by expanding the expression $\left. \frac{\partial \phi}{\partial x}\right|_{j+\frac{3}{2}}$ in Eqn. \ref{eqn:4thOrderOpti}c. The resulting coefficients can be expressed as follows:
\begin{equation} \label{bpcoeffs}
\begin{aligned}
    b_{-3} = \frac{3}{640} - &\psi_{\frac{3}{2}}, \text{ } b_{-2} = \frac{3}{640} (-7 + 1280 \psi_{\frac{3}{2}}), \text{ } b_{-1} = \frac{-3}{32} (-1 + 160 \psi_{\frac{3}{2}}), \text{ } b_{0} = \frac{1}{192} (-19 + 3840 \psi_{\frac{3}{2}}), \\
    &b_{1} = \frac{-3}{128} (43 + 640 \psi_{\frac{3}{2}}), \text{ } b_{2} = \frac{3}{640} (229 + 1280 \psi_{\frac{3}{2}}), \text{ } b_{3} = \frac{1}{960} (-31 - 960 \psi_{\frac{3}{2}})
\end{aligned}
\end{equation}

Substituting the computed midpoint derivative approximations obtained using the aforementioned coefficients $a_p$ and $b_p$ into Eqn. \ref{gen2ndDer}, yields the following expression for the second derivative:
\begin{equation}\label{4thSecondDer}
    \begin{aligned}
        \left.\frac{\partial^2 \phi}{\partial x^2}\right|_{j}^{\text{ME4-Opti}} = 
        & \left( \frac{\phi_{j+1}-2\phi_{j}+\phi_{j-1}}{2\Delta x}\right) \frac{183-160\left(27\psi_{\frac{1}{2}}-\psi_{\frac{3}{2}}\right)}{64}+ \\
        & \left( \frac{\phi_{j+2}-2\phi_{j}+\phi_{j-2}}{4\Delta x}\right) \frac{-39+160\left(27\psi_{\frac{1}{2}}-\psi_{\frac{3}{2}}\right)}{80}+ \\
        & \left( \frac{\phi_{j+3}-2\phi_{j}+\phi_{j-3}}{6\Delta x}\right) \frac{37-480\left(27\psi_{\frac{1}{2}}-\psi_{\frac{3}{2}}\right)}{960} + \mathcal{O}(\Delta x^4)
    \end{aligned}
\end{equation}


\begin{remark}
    Eqn. \ref{4thSecondDer} represents a family of fourth-order second derivative schemes with arbitrary spectral properties. Despite the utilization of partially one-sided midpoint derivatives at the boundaries $j-\frac{3}{2}$ and $j+\frac{3}{2}$ (Eqn. \ref{eqn:4thOrderOpti}), it is noteworthy that the final second derivative scheme presented in Eqn. \ref{4thSecondDer} is entirely central in nature with respect to node $j$ due to the opposing symmetry of the $\frac{\partial \phi}{\partial x}$ stencils at the boundaries.
\end{remark}


Building upon the Fourier error analysis introduced in Sec. \ref{sec:baseStraight},  the following parametric modified wavenumber expression is derived for Eqn. \ref{4thSecondDer}:

\begin{equation}
    \begin{aligned}
        \left(k^{*}_{xx}\right)^{\text{ME4-Opti}}=\frac{1}{720} \Bigr[& -3471+38880 \psi_{\frac{1}{2}}-1440 \psi_{\frac{3}{2}} + \\
        & \left(628-51840 \psi_{\frac{1}{2}}+1920 \psi_{\frac{3}{2}}\right) \cos (k)+ \\
        & \left.\left(-37+12960 \psi_{\frac{1}{2}}-480 \psi_{\frac{3}{2}}\right) \cos (2 k)\right]  \sin ^2\left(\frac{k}{2}\right)
    \end{aligned}
\end{equation}

By varying the parameters $\psi_\frac{1}{2}$ and $\psi_\frac{3}{2}$ in the above expression, a brute-force grid search procedure is employed to maximize the spectral resolving efficiency. This spectral resolving efficiency, denoted as $e_v$, serves as the objective function for the optimization process. It can be defined as the smallest wavenumber ($k$) at which $k^{*}_{xx}$ of the scheme deviates from the exact value ($k^{*,\text{ex}}_{xx}=-k^2$) by a specified error tolerance. The error tolerance is set to five percent in the present study; i.e., $\epsilon = 0.05$. The resolving efficiency, as presented in Ref. \cite{lele1992compact}, can be expressed as follows:

\begin{equation}
    \frac{|k^{*}_{xx}-k^{*,\text{ex}}_{xx}|}{k^{*,\text{ex}}_{xx}} \geq \epsilon
\end{equation}

To further ensure an accurate representation of the computed second derivative in the Fourier error analysis, a constraint is imposed on $k^{*}_{xx}$ during the search process. This constraint ensures that the deviation of local dissipation error (signed error) remains above a predetermined threshold of five percent in relation to the exact theoretical value at all wavenumbers, i.e. $k^{\prime \prime} > 1.05k^{*,\text{ex}}_{xx}, \forall \text{ } k \in [0,\pi]$.\\

\begin{figure}[h!]
    \centering
    \includegraphics[width=80mm]{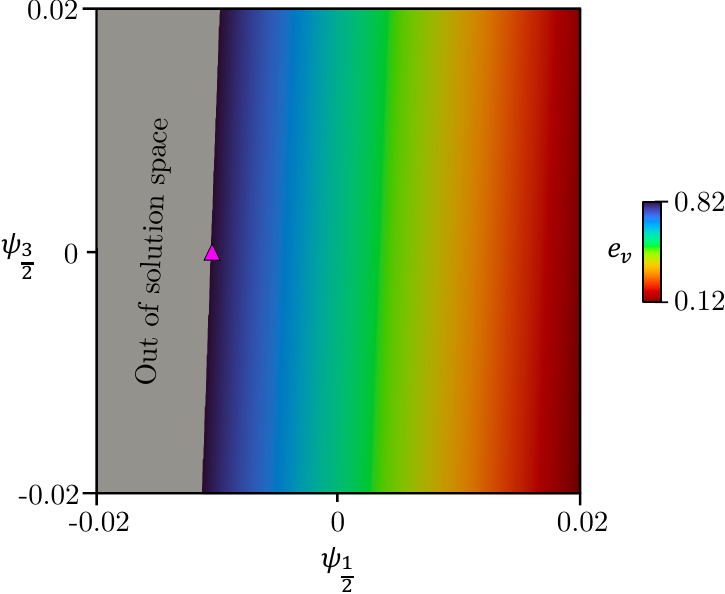}
    \caption{Contour distribution of resolving efficiency ($e_v$) of the second derivative discretization, plotted as a function of leading truncation error term coefficients $\psi_\frac{1}{2}$ and $\psi_\frac{3}{2}$ that appear in Eqn. \ref{eqn:4thOrderOpti}. The magenta triangle indicates the optimal solution at $\psi_\frac{1}{2}=-0.01$ and $\psi_\frac{3}{2}=0$, where the resolving efficiency reaches its maximum.}
    \label{OptiPlot}
\end{figure}

Fig. \ref{OptiPlot} illustrates the search space employed to find the optimal values of the coefficients $\psi_\frac{1}{2}$ and $\psi_\frac{3}{2}$, along with the corresponding resolving efficiency ($e_v$) of the scheme within this parameter space. The resolving efficiency demonstrates a higher sensitivity to variations in $\psi_\frac{1}{2}$ compared to $\psi_\frac{3}{2}$. Using the optimal solution $\psi_\frac{1}{2}=-0.01$ and $\psi_\frac{3}{2}=0$ obtained from Fig. \ref{OptiPlot}, the coefficients $a_p$ and $b_p$ are computed using Eqns. \ref{apcoeffs} and \ref{bpcoeffs} respectively. The resulting values of $a_p$ and $b_p$ are summarized in Table \ref{table:4thOrderopti}.

\begin{table}[h!]
    \centering
    \caption{Discretization coefficients for the Explicit 4th order Optimized scheme ME4-Opti.}
    \begin{tabular}{cccccccc}
    \toprule & $p=-3$ & $p=-2$ & $p=-1$ & $p=0$ & $p=1$ & $p=2$ & $p=3$ \\
    \midrule
    $a_p$ & $\frac{133}{12500}$ & $\frac{-27411}{400000}$ & $\frac{53929}{240000}$ & $\frac{-55387}{40000}$ & $\frac{53259}{40000}$ & $\frac{-154733}{1200000}$ & $\frac{6131}{400000}$ \\
    $b_p$ & $\frac{623}{80000}$ & $\frac{-4113}{80000}$ & $\frac{561}{4000}$ & $\frac{-3863}{24000}$ & $\frac{-15381}{16000}$ & $\frac{84387}{80000}$ & $\frac{-3503}{120000}$ \\
    \bottomrule
    \end{tabular}
    \label{table:4thOrderopti}
\end{table}

This concludes the presentation of the fourth order optimized straight second derivative scheme. The procedure described so far to derive the fourth order scheme, has been followed again to obtain the sixth order optimized second derivative scheme ME6-Opti. The seventh order nine point long stencil expressions used to evaluate the midpoint derivatives for this scheme are presented in the following Eqn. \ref{eqn:6thOrderOpti}:

\begin{equation}\label{eqn:6thOrderOpti}
    \begin{aligned}
    \left. \frac{\partial \phi}{\partial x}\right|_{j \pm \frac{1}{2}} &= \sum_{p=-4}^{4} a_{ \pm p} \frac{\pm\phi_{j+p}}{\Delta x}  \pm  \vartheta_{\frac{1}{2}} \frac{\partial^7\phi}{\partial x^7} \Delta x^7 + \mathcal{O}(\Delta x^8), \\
    \left. \frac{\partial \phi}{\partial x}\right|_{j \pm \frac{3}{2}} & = \sum_{p=-4}^{4} b_{ \pm p} \frac{\pm\phi_{j+p}}{\Delta x}  \pm  \vartheta_{\frac{3}{2}} \frac{\partial^7\phi}{\partial x^7} \Delta x^7 + \mathcal{O}(\Delta x^8), \\
    \left. \frac{\partial \phi}{\partial x}\right|_{j \pm \frac{5}{2}} &= \sum_{p=-4}^{4} c_{ \pm p} \frac{\pm\phi_{j+p}}{\Delta x}  \pm  \vartheta_{\frac{5}{2}} \frac{\partial^7\phi}{\partial x^7} \Delta x^7 + \mathcal{O}(\Delta x^8). \\
    \end{aligned}
    \left\}
    \begin{aligned}
    \\
    \\
    \\
    \\
    \\
    \\
    \\
    \end{aligned}
    \right.
    \text{7th order}
\end{equation}

Where $\vartheta_{\frac{1}{2}}$, $\vartheta_{\frac{3}{2}}$, and $\vartheta_{\frac{5}{2}}$ are the coefficients of leading seventh order truncation error terms in the Taylor series expansion. Substituting these midpoint derivative expressions (Eqn. \ref{eqn:6thOrderOpti}) into Eqn. \ref{gen2ndDer}, will yield the following family of sixth order central second derivative expression which is subject to optimization under the same constraints as that of the fourth order scheme: 
\begin{equation}\label{6thsecondDer}
    \begin{aligned}
        \left.\frac{\partial^2 \phi}{\partial x^2}\right|_{j}^{\text{ME6-Opti}} = 
        & \left( \frac{\phi_{j+1}-2\phi_{j}+\phi_{j-1}}{2\Delta x}\right) \frac{2997+112\left(2250\vartheta_{\frac{1}{2}}-125\vartheta_{\frac{3}{2}}+9\vartheta_{\frac{5}{2}}\right)}{960}+ \\
        & \left( \frac{\phi_{j+2}-2\phi_{j}+\phi_{j-2}}{4\Delta x}\right) \frac{-693-112\left(2250\vartheta_{\frac{1}{2}}-125\vartheta_{\frac{3}{2}}+9\vartheta_{\frac{5}{2}}\right)}{960}+ \\
        & \left( \frac{\phi_{j+3}-2\phi_{j}+\phi_{j-3}}{6\Delta x}\right) \frac{799+336\left(2250\vartheta_{\frac{1}{2}}-125\vartheta_{\frac{3}{2}}+9\vartheta_{\frac{5}{2}}\right)}{6720}+ \\
        & \left( \frac{\phi_{j+4}-2\phi_{j}+\phi_{j-4}}{8\Delta x}\right) \frac{-117-112\left(2250\vartheta_{\frac{1}{2}}-125\vartheta_{\frac{3}{2}}+9\vartheta_{\frac{5}{2}}\right)}{5760} + \mathcal{O}(\Delta x^6)
    \end{aligned}
\end{equation}

After the optimization process, the leading term truncation error coefficients $\vartheta_{\frac{1}{2}}$, $\vartheta_{\frac{3}{2}}$, and $\vartheta_{\frac{5}{2}}$ are determined to be $3/1250$, $-1/1250$, and $3/625$ respectively. The subsequent optimized coefficients $a_p$, $b_p$, and $c_p$ are summarized in Table \ref{6thOptiTable}.




\begin{table}[h!]
    \centering
    \caption{Discretization coefficients for the Explicit 6th order Optimized scheme ME6-Opti.}
    \begin{tabular}{cccccccccc}
    \toprule & $p=-4$ & $p=-3$ & $p=-2$ & $p=-1$ & $p=0$ & $p=1$ & $p=2$ & $p=3$ & $p=4$ \\
    \midrule
    $a_p$ & $\frac{-3}{1250}$ & $\frac{89141}{4480000}$ & $\frac{-49133}{640000}$ & $\frac{411173}{1920000}$ & $\frac{-174629}{128000}$ & $\frac{851641}{640000}$ & $\frac{-282149}{1920000}$ & $\frac{18413}{640000}$ & $\frac{-13877}{4480000}$ \\
    $b_p$ & $\frac{459}{4480000}$ & $\frac{-547}{4480000}$ & $\frac{-1289}{640000}$ & $\frac{2703}{640000}$ & $\frac{18379}{384000}$ & $\frac{-738047}{640000}$ & $\frac{742461}{640000}$ & $\frac{-820391}{13440000}$ & $\frac{9167}{2240000}$ \\
    $c_p$ & $\frac{-3377}{2240000}$ & $\frac{36157}{4480000}$ & $\frac{-6141}{640000}$ & $\frac{-20593}{640000}$ & $\frac{16367}{128000}$ & $\frac{-296029}{1920000}$ & $\frac{-618391}{640000}$ & $\frac{4737907}{4480000}$ & $\frac{-4000637}{13440000}$ \\
    \bottomrule
    \end{tabular}
    \label{6thOptiTable}
\end{table}

By utilizing the optimized coefficients $a_p$, $b_p$, and $c_p$ provided in Tables \ref{table:4thOrderopti} and \ref{6thOptiTable}, the Fourier analysis of the ME4-Opti and ME6-Opti second derivative schemes yields the following modified wavenumbers ($k^{*}_{xx}$) which dictate their respective spectral properties.

\begin{equation}
    \begin{aligned}
        \left(k^{*}_{xx}\right)^{\text{ME4-Opti}} &= -k^2 \left(1-\frac{2437 k^4}{240000}+\frac{5059 k^6}{6720000}-\frac{353 k^{8}}{19200000}+ \cdots\right) \\
        \left(k^{*}_{xx}\right)^{\text{ME6-Opti}} &= -k^2 \left(1+\frac{341149 k^6}{67200000}-\frac{405149 k^{8}}{201600000}+\frac{90500141 k^{10}}{266112000000}- \cdots \right) 
    \end{aligned}
\end{equation} \label{ME46Opti-modified}

\begin{figure}[h!]
    \centering
    \includegraphics[width=150mm]{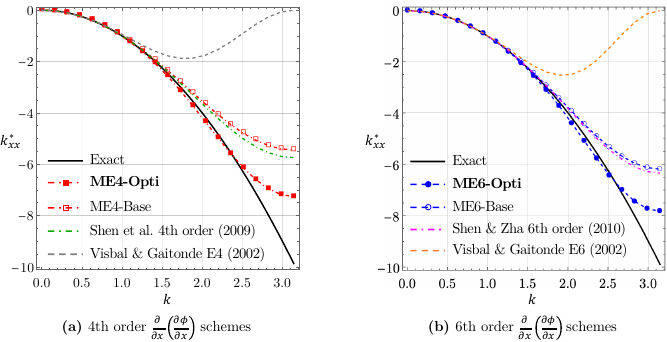}
    \caption{Modified wavenumber profiles of optimized 4th and 6th order straight derivative term discretization schemes, in comparison to the baseline and reference schemes from the literature of same accuracy order.}
    \label{optimizedSchemes}
\end{figure}

Fig. \ref{optimizedSchemes}(a,b) presents a comparison of the spectral properties between the optimized schemes and the corresponding baseline and reference schemes of the same order found in the literature. It is evident that the optimized schemes exhibit a significant improvement in their spectral properties when compared to the baseline schemes. Also, it is noteworthy that the modified wavenumber profiles for the Shen et al. 4th order scheme \cite{shen2009high} and Shen \& Zha 6th order scheme \cite{shen2010large} closely resemble the profiles of the baseline schemes.

To assess the dissipation accuracy of the derived schemes, the spectral viscosity associated with each scheme is quantified at the Nyquist limit. The spectral viscosity defined as $\nu^{\prime \prime}_s = -\nu (k^{*}_{xx}+k^2)/k^2$, is the effective artificial numerical viscosity introduced in to the system of governing equations as a result of the numerical error at a given wavenumber. While an ideal scenario would have a spectral viscosity of zero at all the wavenumbers, in practice, it can take either positive or negative value (anti-dissipation) depending on the local wavenumber and the nature of the scheme employed. Lambalias et al. \cite{lamballais2011straightforward, lamballais2021viscous} have demonstrated the customization of spectral viscosity profiles to introduce targeted numerical dissipation in the high wavenumber range, addressing sub-grid scale dissipation or creating an up-winding effect. Nonetheless, the focus of the schemes developed in this paper is primarily on accurately representing the second derivative across the entire wavenumber range while maintaining the desired order of accuracy. \textcolor{magenta}{The equivalent numerical dissipation and corresponding equivalent Reynolds number due to discretization error of second derivative schemes are discussed in detail in \ref{app:eqivalRe}.}\\

\begin{table}[h!] 
    \centering
    \caption{Resolving efficiency $e_v$ under $5\%$ error tolerance limit (i.e. $\epsilon=0.05$) and spectral viscosity $\nu_s^{\prime \prime}$ at cutoff wavenumber for baseline, optimized, and reference second derivative discretization schemes from the literature.}
    \includegraphics[width=135mm]{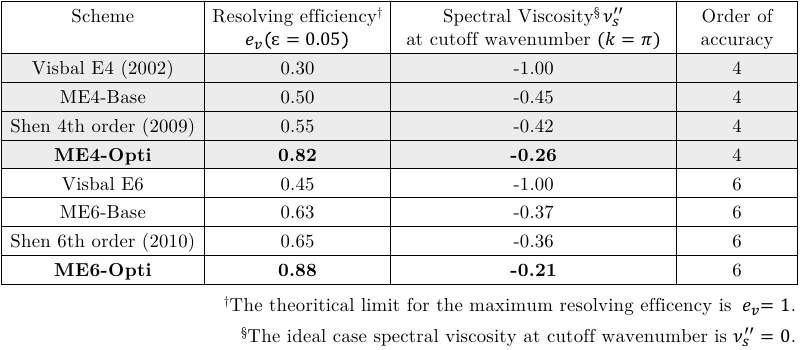}
    \label{ResolvingEff}
\end{table}

Table \ref{ResolvingEff} provides a summary of the resolving efficiency $e_v$ and spectral viscosity $\nu^{\prime \prime}_s$ at the cutoff wavenumber $k=\pi$ (Nyquist limit) for various schemes. It is noteworthy that both the optimized schemes derived in this section exhibit high resolving efficiencies, exceeding 80\%. These optimized schemes demonstrate a modest resolving efficiency improvement of $\approx$32\% and $\approx$25\% for the 4th and 6th order schemes respectively, despite employing a smaller stencil for the 6th order scheme ME6-Opti. This is attributed to the efficient use of information within the stencil to compute derivatives at each midpoint locations and the optimization procedure. While the ME4-Opti and ME6-Opti schemes exhibit minimal magnitudes of spectral viscosity of all the schemes compared (thus most accurate), the Visbal E4/E6 schemes display a spectral viscosity of $\nu^{\prime \prime}=-1.0$ at the cutoff wavenumber. This is due to the absence of the h-elliptic property \cite{nishikawa2010beyond} and the resulting odd-even decoupling effect. In addition to its impact on accurately diffusing the high-wavenumber features, the h-elliptic property has been identified in previous studies to also enhance the numerical stability. The impact of h-elliptic property and spectral properties in particular on the computed flow solutions is detailed in \cite{chamarthi2022importance,chamarthi2023role}. In summary, the derived 4th and 6th order optimized schemes are \textit{high resolution} and \textit{odd-even mode suppressing} in nature. Furthermore, the ME6-Opti scheme stands out for its utilization of a smaller stencil, making it \textit{compact} and \textit{efficient}.

\subsection{Discretization of mixed derivative terms} \label{sec:mixedDers}
Thus far, the discretization of straight derivative terms (e.g., of the form $\frac{\partial}{\partial x} \cdot \frac{\partial \phi}{\partial x}$) has been discussed. However, in Eqn. \ref{tauxxTerm}, the presence of tangential shear stress and bulk viscosity terms $-\mu\frac{2}{3}\frac{\partial v}{\partial y}$ and $-\mu \frac{2}{3}\frac{\partial w}{\partial z}$ will introduce mixed second derivative terms into the calculation. The baseline and optimized ME4/ME6 discretization schemes for mixed derivative terms are presented now. The optimized schemes feature high-order accuracy and a `filter penalty term' for accurate spectral representation in the high wavenumber range. \\

Dropping the viscous coefficient, the following expression can be used to generalize various mixed derivative terms present in the governing equations:


\begin{equation}
    \frac{\partial g}{\partial x} = \frac{\partial }{\partial x} \left(\frac{\partial \phi}{\partial y} \right)
\end{equation}

Similar to Eqn \ref{nonlinear1}, the outer derivative $\frac{\partial}{\partial x}$ in the above expression can be discretized as follows:

\begin{equation} \label{mixedDer}
    \begin{aligned}
        \left. \frac{\partial g}{\partial x} \right|_{j,l} \approx \frac{a^*}{\Delta x} & \left[ \left( \frac{\partial \phi}{\partial y}\right)_{j+\frac{1}{2},l} - \left( \frac{\partial \phi}{\partial y}\right)_{j-\frac{1}{2},l}  \right] + \frac{b^*}{3 \Delta x} \left[ \left( \frac{\partial \phi}{\partial y}\right)_{j+\frac{3}{2},l} - \left( \frac{\partial \phi}{\partial y}\right)_{j-\frac{3}{2},l}  \right] + \\
        \frac{c^*}{5 \Delta x} & \left[ \left( \frac{\partial \phi}{\partial y}\right)_{j+\frac{5}{2},l} - \left( \frac{\partial \phi}{\partial y}\right)_{j-\frac{5}{2},l}  \right],
    \end{aligned}
\end{equation}

\noindent where $l$ refers to the index along the $y$-direction. To compute terms on the right-hand side (RHS) of the above Eqn. \ref{mixedDer}, we begin by approximating the nodal values $\left . \frac{\partial \phi}{\partial y} \right|_{j,l}$ using either central fourth-order differences for the ME4 scheme or central sixth-order differences for the ME6 scheme:

\begin{equation} \label{dfdyDer}
    \begin{aligned}
        \left. \frac{\partial \phi}{\partial y} \right|_{j,l} &= \frac{2}{3 \Delta y} \left( \phi_{j,l+1} - \phi_{j,l-1} \right) - \frac{1}{12 \Delta y} \left( \phi_{j,l+2} - \phi_{j,l-2} \right) + \mathcal{O}(\Delta y^4)\\
        \left. \frac{\partial \phi}{\partial y} \right|_{j,l} &= \frac{3}{4 \Delta y} \left( \phi_{j,l+1} - \phi_{j,l-1} \right) - \frac{3}{20 \Delta y} \left( \phi_{j,l+2} - \phi_{j,l-2} \right) + \frac{1}{60 \Delta y} \left( \phi_{j,l+3} - \phi_{j,l-3} \right) + \mathcal{O}(\Delta y^6)
    \end{aligned}
\end{equation}

The cell center derivatives $\left . \frac{\partial \phi}{\partial y} \right|_{j,l}$ are interpolated to half locations using two different approaches. The first approach corresponds to the baseline (ME4-Base and ME6-Base), while the second approach corresponds to the optimized method (ME4-Opti and ME6-Opti). These two are detailed in the following two subsections \ref{sec:mixed1} and \ref{sec:mixed2}.

\subsubsection{ME4 and ME6 baseline schemes for mixed derivative terms} \label{sec:mixed1}

In the baseline ME4 and ME6 schemes, the values of $\frac{\partial \phi}{\partial y}$ computed at each $(j,l)$ are directly interpolated to the midpoint locations using the following standard fourth-order (for ME4-Base) or sixth-order (for ME6-Base) interpolation formulas:

\begin{equation}\label{InterpolYder-base}
    \begin{aligned}
        \left. \frac{\partial \phi}{\partial y} \right|_{j+\frac{1}{2}} &= \frac{9}{16} \left[ \left( \frac{\partial \phi}{\partial y}\right)_{j+1} + \left( \frac{\partial \phi}{\partial y}\right)_{j} \right] - \frac{1}{16} \left[ \left( \frac{\partial \phi}{\partial y}\right)_{j+2} + \left( \frac{\partial \phi}{\partial y}\right)_{j-1} \right] + \mathcal{O}\left( \Delta x^4 \right) \\
        \left. \frac{\partial \phi}{\partial y} \right|_{j+\frac{1}{2}} &= \frac{75}{128} \left[ \left( \frac{\partial \phi}{\partial y}\right)_{j+1} + \left( \frac{\partial \phi}{\partial y}\right)_{j} \right] - \frac{25}{256} \left[ \left( \frac{\partial \phi}{\partial y}\right)_{j+2} + \left( \frac{\partial \phi}{\partial y}\right)_{j-1} \right] \\
        & + \frac{3}{256} \left[ \left( \frac{\partial \phi}{\partial y}\right)_{j+3} + \left( \frac{\partial \phi}{\partial y}\right)_{j-2} \right] + \mathcal{O}\left( \Delta x^6 \right)
    \end{aligned}
\end{equation}

Interpolation to other midpoint locations e.g., $\left[j+\frac{3}{2}, j+\frac{5}{2}, \ldots \right]$ can be performed using the same formulae by shifting the stencil. The resultant interpolated derivative values are subsequently substituted in Eqn. \ref{mixedDer} and the mixed derivative is evaluated.

To evaluate the spectral representation of the derived mixed derivative schemes, a 2-D Fourier error analysis is conducted. As these mixed second derivative terms involve both $x$ and $y$ spatial directions, we consider a two-dimensional Fourier basis function $\mathcal{F}(x,y)$ with wavenumber $k$, represented by the following Eqn. \ref{2Dfourer}.

\begin{equation} \label{2Dfourer}
    \mathcal{F}(x,y)=e^{i k (x+y)}, \quad k=\frac{2 \pi}{L} n, \quad n=0,1,2, \ldots, N / 2
\end{equation}

The error analysis is conducted on a uniformly spaced 2-D domain of length `$L$' with a grid containing $N$ points in each direction. The analytical expression for the mixed second derivative is as follows:

\begin{equation}
    \frac{\partial }{\partial x} \left( \frac{\partial \mathcal{F}(x,y)}{\partial y} \right) = \mathcal{F}_{xy}(x,y) = -k^2 \mathcal{}
\end{equation}

\noindent The modified wavenumber $k^{*}_{xy}$ is defined as follows:

\begin{equation}
    k^{*}_{xy} = \frac{\mathcal{F}_{xy}(x,y)}{\mathcal{F}(x,y)} = -k^2
\end{equation}

Similarly, the discrete values of $\mathcal{F}(x,y)$ on the 2-D grid are used to derive modified wavenumber expressions for both the ME4-Base and ME6-Base mixed derivative schemes. This process entails substituting the values of Fourier modes $\mathcal{F}(x,y)$ associated with different grid nodes into the respective schemes, thereby obtaining the following modified wavenumber expressions for ME4-Base and ME6-Base schemes respectively:

\begin{equation}
    \begin{aligned}
        \left(k^{*}_{xy}\right)^{\text{ME4-Base}} &= -\frac{59}{64} + \frac{275 \cos (k)}{1152}+\frac{65 \cos (2 k)}{72} -\frac{61 \cos (3 k)}{256} +\frac{11 \cos (4 k)}{576} -\frac{\cos (5 k)}{2304}\\
        &=-k^2 \left( 1 - \frac{59 k^4}{960} + \frac{25 k^6}{4032} + \frac{503 k^8}{691200} + \cdots \right)\\
    \end{aligned}
\end{equation}

\begin{equation}
    \begin{aligned}
        \left(k^{*}_{xy}\right)^{\text{ME6-Base}} &=  -\frac{704663}{589824} + \frac{31895 \cos (k)}{65536}+\frac{333251 \cos (2 k)}{307200}-\frac{774123 \cos (3 k)}{1638400}+\frac{26711 \cos (4 k)}{245760} \\
        & \quad -\frac{2779 \cos (5 k)}{196608}+\frac{223 \cos (6 k)}{184320}-\frac{281 \cos(7 k)}{4915200}+\frac{3 \cos (8 k)}{1638400}\\
        &=-k^2 \left(1 - \frac{57 k^6}{4480} + \frac{1649 k^8}{737280} - \frac{5309 k^{10}}{27033600} + \cdots \right)
    \end{aligned}
\end{equation}

\begin{remark}
    It can be observed that the central nature of the derivative formulae in Eqn. \ref{dfdyDer} for computing $\frac{\partial \phi}{\partial y}$ results in the final expression for the mixed derivative term $\left. \frac{\partial g}{\partial x} \right|_{j,l}$ in Eqn. \ref{mixedDer}, which does not contain any information belonging to the index $l$. Instead, only the indices $l \pm 1$, $l \pm 2$, and $l \pm 3$ are present in the final formula. This limitation leads to poor mixed derivative spectral properties. In the optimized mixed derivative schemes, we resolve this issue by introducing a \textit{filter penalty term} that includes points belonging to the index $l$. The magnitude of damping provided by this filter term is controlled by adjusting the leading error terms in the Taylor expansion of the filter expression.
\end{remark}

\subsubsection{ME4 and ME6 optimized schemes for mixed derivative terms} \label{sec:mixed2}
In the present approach, along with an interpolation term, an additional filter penalty term is employed in the formulation unlike the base scheme which only contains an interpolation term. The general interpolation expression for both ME4-Opti and ME6-Opti is written as follows:


\begin{equation}\label{eqn:InterpFilter}
    \left. \frac{\partial \phi}{\partial y}\right|_{j \pm \frac{1}{2},l} = \underbrace{\sum_{p=-n}^{n} a_{ \pm p}^{I} \left(\frac{\partial \phi}{\partial y}\right)_{j \pm p,l}}_{\text{Interpolation term}} + \underbrace{\sum_{p=-n}^{n} \frac{a_{ \pm p}^{\varnothing}}{\Delta x} \phi_{j + p,l}}_{\text{Filter penalty term}} + \mathcal{O}(\Delta x^{2n-2}, \Delta y^{2n-2})
\end{equation}

Similarly, for $j \pm \frac{3}{2}$ and $j \pm \frac{5}{2}$ interpolation, the coefficient set $(a_{ \pm p}^{I}, a_{ \pm p}^{\varnothing})$ is replaced respectively with $(b_{ \pm p}^{I}, b_{ \pm p}^{\varnothing})$ and $(c_{ \pm p}^{I}, c_{ \pm p}^{\varnothing})$ in the above expression. The values of these coefficients are obtained through the same optimization procedure that was previously discussed in the context of straight derivative terms in Sec. \ref{sec:optiStraight} and is not repeated here. The obtained coefficients for fourth and sixth order are summarized in tables \ref{mixedTable4th} and \ref{mixedTable6th}.



\begin{table}[h!]
    \centering
    \begin{tabular}{cccccccc}
        \toprule
        & $p=-3$ & $p=-2$ & $p=-1$ & $p=0$ & $p=1$ & $p=2$ & $p=3$ \\
        \midrule
        $a_p^I$ & $\frac{-83}{384000}$ & $\frac{1473}{64000}$ & $\frac{-21363}{128000}$ & $\frac{72409}{96000}$ & $\frac{49497}{128000}$ & $\frac{1129}{64000}$ & $\frac{-5567}{384000}$ \\
        $b_p^I$ & $\frac{811}{128000}$ & $\frac{-3151}{64000}$ & $\frac{4469}{25600}$ & $\frac{-2529}{6400}$ & $\frac{4661}{5120}$ & $\frac{23977}{64000}$ & $\frac{-2753}{128000}$ \\
        $a_p^{\varnothing}$ & $\frac{1}{20}$ & $-\frac{3}{10}$ & $\frac{3}{4}$ & -1 & $\frac{3}{4}$ & $-\frac{3}{10}$ & $\frac{1}{20}$ \\
        $b_p^{\varnothing}$ & $\frac{-1}{2000}$ & $\frac{3}{1000}$ & $\frac{-3}{400}$ & $\frac{1}{100}$ & $\frac{-3}{400}$ & $\frac{3}{1000}$ & $\frac{-1}{2000}$ \\
        \bottomrule
    \end{tabular}
    \caption{Coefficients of the fourth-order mixed derivative optimized scheme (ME4-Opti) with $n=3$, as represented by Eqn. \ref{eqn:InterpFilter}.}
    \label{mixedTable4th}
\end{table}

\begin{table}[h!]
    \centering
    \begin{tabular}{cccccccccc}
    \toprule
    & $p=-4$ & $p=-3$ & $p=-2$ & $p=-1$ & $p=0$ & $p=1$ & $p=2$ & $p=3$ & $p=4$ \\
    \midrule
    $a_p^I$ & $\frac{-661}{819200}$ & $\frac{263}{512000}$ & $\frac{31573}{1024000}$ & $\frac{-91107}{512000}$ & $\frac{302761}{409600}$ & $\frac{43093}{102400}$ & $\frac{6429}{1024000}$ & $\frac{-12349}{512000}$ & $\frac{21511}{4096000}$ \\
    $b_p^I$ & $\frac{-7673}{4096000}$ & $\frac{9179}{512000}$ & $\frac{-15959}{204800}$ & $\frac{106337}{512000}$ & $\frac{-165879}{409600}$ & $\frac{456421}{512000}$ & $\frac{408037}{1024000}$ & $\frac{-3357}{102400}$ & $\frac{8279}{4096000}$ \\
    $c_p^I$ & $\frac{8279}{4096000}$ & $\frac{-10273}{512000}$ & $\frac{92869}{1024000}$ & $\frac{-126827}{512000}$ & $\frac{37877}{81920}$ & $\frac{-337743}{512000}$ & $\frac{1086701}{1024000}$ & $\frac{166763}{512000}$ & $\frac{-59769}{4096000}$ \\
    $a_p^{\varnothing}$ & $-\frac{13}{1000}$ & $\frac{13}{125}$ & $-\frac{91}{250}$ & $\frac{91}{125}$ & $-\frac{91}{100}$ & $\frac{91}{125}$ & $-\frac{91}{250}$ & $\frac{13}{125}$ & $-\frac{13}{1000}$ \\
    $b_p^{\varnothing}$ & $-\left(\frac{1}{2000}\right)$ & $\frac{1}{250}$ & $-\left(\frac{7}{500}\right)$ & $\frac{7}{250}$ & $-\left(\frac{7}{200}\right)$ & $\frac{7}{250}$ & $-\left(\frac{7}{500}\right)$ & $\frac{1}{250}$ & $-\left(\frac{1}{2000}\right)$ \\
    $c_p^{\varnothing}$ & $-\left(\frac{1}{2000}\right)$ & $\frac{1}{250}$ & $-\left(\frac{7}{500}\right)$ & $\frac{7}{250}$ & $-\left(\frac{7}{200}\right)$ & $\frac{7}{250}$ & $-\left(\frac{7}{500}\right)$ & $\frac{1}{250}$ & $-\left(\frac{1}{2000}\right)$ \\
    \bottomrule
    \end{tabular}
    \caption{Coefficients of the sixth-order mixed derivative optimized scheme (ME6-Opti) with $n=4$, as represented by Eqn. \ref{eqn:InterpFilter}.}
    \label{mixedTable6th}
\end{table}

Using the listed coefficients in tables \ref{mixedTable4th} and \ref{mixedTable6th}, the Fourier error analysis of the ME4-Opti and ME6-Opti is conducted. The resultant expressions for the modified wavenumber $k_{xy}^{*}$ are as follows:

\begin{equation}
    \begin{aligned}
        \left(k^{*}_{xy}\right)^{\text{ME4-Opti}} &= -\frac{438379}{144000}+\frac{1009171 \cos (k)}{288000}-\frac{487 \cos (2 k)}{900} +\frac{10919 \cos (3 k)}{115200}-\frac{2293 \cos (4 k)}{144000}+\frac{159 \cos (5 k)}{64000}\\
        &=-k^2 \left( 1 - \frac{257 k^4}{16000} + \frac{20687 k^6}{2016000} - \frac{120101 k^{8}}{34560000} + \cdots \right)
    \end{aligned}
\end{equation}

\begin{equation}
    \begin{aligned}
        \left(k^{*}_{xy}\right)^{\text{ME6-Opti}} &=  -\frac{180127829}{57600000}+ \frac{28259327 \cos (k)}{7680000}-\frac{81089207 \cos (2 k)}{115200000}+\frac{7562747 \cos (3 k)}{38400000}-\frac{671839 \cos (4 k)}{11520000} \\
        & \quad +\frac{1784983 \cos (5 k)}{115200000}-\frac{65173 \cos (6k)}{23040000}+\frac{25991 \cos (7 k)}{115200000}\\
        &=-k^2 \left(1 + \frac{34751 k^6}{13440000} - \frac{39581 k^8}{23040000} + \frac{815231 k^{10}}{3041280000} + \cdots \right)
    \end{aligned}
\end{equation}

\begin{figure}[h!]
    \centering
    \includegraphics[width=150mm]{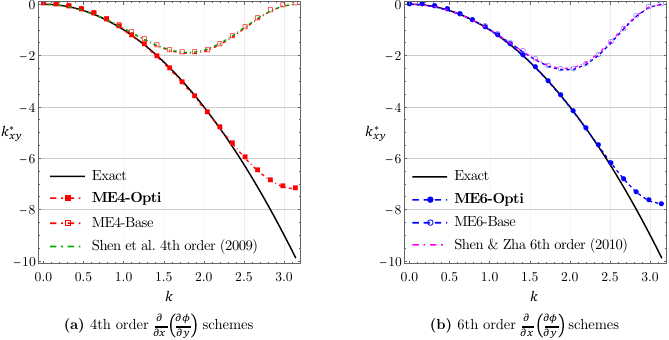}
    \caption{Modified wavenumber profiles of the optimized 4th and 6th order mixed derivative discretization schemes compared with the baseline and reference schemes from literature.}
    \label{MixedDerSpectral}
\end{figure}

The leading term in the expansion series presented above demonstrates the formal fourth and sixth-order accuracy of the respective schemes. The plots of the modified wavenumber expressions obtained are shown in Fig. \ref{MixedDerSpectral}. The superior spectral representation in the wavenumber range $k=[\pi/2,\pi]$ of the present optimized schemes in comparison to the baseline and the schemes from literature can be clearly appreciated. This is attributed to the filtering term added in to the Eqn. \ref{eqn:InterpFilter}. On the other hand, the baseline and the schemes of Shen et al. \cite{shen2009high} and Shen \& Zha \cite{shen2010large} resemble each other very closely and do not possess good spectral representation at high wavenumbers.

\subsection{Viscosity coefficient at midpoint locations} \label{sec:viscCoeff}
The past two sections \ref{sec:straight} and \ref{sec:mixedDers} have discussed the discretization of the linear second derivative terms of form $\frac{\partial }{\partial x} \left(\frac{\partial \phi}{\partial x} \right)$ (straight) and $\frac{\partial }{\partial x} \left(\frac{\partial \phi}{\partial y} \right)$ (mixed). However, the viscous flux terms in the Navier-Stokes equations are indeed non-linear in nature. For instance the term $\frac{\partial }{\partial x} \left( \frac{4 \mu}{3} \frac{\partial u}{\partial x} \right)$ in Eqn. \ref{tauxxTerm} is non-linear under variable viscosity conditions. Such expressions can be generalized using the term $\frac{\partial }{\partial x} \left( \mu \frac{\partial \phi}{\partial x} \right)$. Just as Eqn. \ref{gen2ndDer} and Eqn. \ref{mixedDer}, these type of non-linear terms can also discretized as follows:



\begin{equation}\label{nonlinear1}
    \begin{aligned}
        \left. \frac{\partial f}{\partial x} \right|_{j} = &\frac{a^*}{\Delta x} \left[ \mu_{j+\frac{1}{2}} \left( \frac{\partial \phi}{\partial x}\right)_{j+\frac{1}{2}} - \mu_{j-\frac{1}{2}} \left( \frac{\partial \phi}{\partial x}\right)_{j-\frac{1}{2}}  \right] + \frac{b^*}{3 \Delta x} \left[ \mu_{j+\frac{3}{2}} \left( \frac{\partial \phi}{\partial x}\right)_{j+\frac{3}{2}} - \mu_{j-\frac{3}{2}} \left( \frac{\partial \phi}{\partial x}\right)_{j-\frac{3}{2}}  \right] + \\
        &\frac{c^*}{5 \Delta x} \left[ \mu_{j+\frac{5}{2}} \left( \frac{\partial \phi}{\partial x}\right)_{j+\frac{5}{2}} - \mu_{j-\frac{5}{2}} \left( \frac{\partial \phi}{\partial x}\right)_{j-\frac{5}{2}}  \right] + \mathcal{O}(\Delta x^{\text{OA}})
    \end{aligned}
\end{equation}

The computation of the derivative term $\frac{\partial \phi}{\partial x}$ in the above equation at various midpoint locations $[j \pm \frac{1}{2}, j \pm \frac{3}{2}, j \pm \frac{5}{2}]$ has already been discussed in Sections \ref{sec:baseStraight} (baseline) and \ref{sec:optiStraight} (optimized). To maintain high-order accuracy while discretizing Eqn. \ref{nonlinear1}, the values of $\mu$ at the midpoint locations need to be computed with a minimum accuracy of fourth and sixth order to satisfy the accuracy requirements of the ME4 and ME6 base/optimized schemes, respectively. For the baseline schemes ME4-Base and ME6-Base, the following standard central interpolation approximations \cite{gaitonde1998high} are employed:


\begin{equation}\label{InterPMu-base}
    \begin{aligned}
        \left( \mu_{j+\frac{1}{2}} \right)^{\text{ME4-Base}} &= \frac{9}{16} \left( \mu_{j+1} + \mu_{j} \right) - \frac{1}{16} \left( \mu_{j+2} + \mu_{j-1} \right) + \mathcal{O}\left( \Delta x^4 \right) \\
        \left( \mu_{j+\frac{1}{2}} \right)^{\text{ME6-Base}} &= \frac{75}{128} \left( \mu_{j+1} + \mu_{j} \right) - \frac{25}{256} \left( \mu_{j+2} + \mu_{j-1} \right) + \frac{3}{256} \left( \mu_{j+3} + \mu_{j-2} \right) + \mathcal{O}\left( \Delta x^6 \right)
    \end{aligned}
\end{equation}

In contrast to the baseline schemes, the ME4-Opti and ME6-Opti schemes take advantage of all the values of $\mu$ within the stencil to construct the interpolation polynomial. The resultant expressions to compute the midpoint dynamic viscosity values are as follows:

\begin{equation} \label{InterPMu}
        \mu_{j\pm\frac{1}{2}} = \sum_{p=-n}^{n} a_{\pm p}^{I} \mu_{p} + \mathcal{O}\left( \Delta x^{2q-2} \right), \quad 
        \mu_{j\pm\frac{3}{2}} = \sum_{p=-n}^{n} b_{\pm p}^{I} \mu_{p} + \mathcal{O}\left( \Delta x^{2q-2} \right), \quad
        \mu_{j\pm\frac{5}{2}} = \sum_{p=-n}^{n} c_{\pm p}^{I} \mu_{p} + \mathcal{O}\left( \Delta x^{2q-2} \right)
\end{equation}

The stencil width is described by $n$, with $n=3$ for the ME4-Opti scheme and $n=4$ for the ME6-Opti scheme. The interpolation coefficients $[a_{\pm p}^{I}, b_{\pm p}^{I}, c_{\pm p}^{I}]$ for the fourth and sixth order optimized schemes are specified in tables \ref{mixedTable4th} and \ref{mixedTable6th} respectively. Analogously, the viscosity is also interpolated via the same formulae Eqn. \ref{InterPMu-base} and \ref{InterPMu} for the case of mixed derivative non-linear terms of form $\frac{\partial }{\partial x} \left( \mu \frac{\partial \phi}{\partial y}\right)$. The boundary conditions are implemented using the ghost cell method avoiding the need of one-sided differencing. More details regarding ghost cell implementation can be found in Ref.\cite{chandravamsi2023application}.






%% file: timeDisc.tex
\section{Time integration and stability criterion} \label{sec:time-int}

\subsection{Time integration} \label{sec:RK3-TI}
The third order explicit three-stage Total Variation Diminishing Runge-Kutta (TVD RK-3) scheme \cite{gottlieb1998total} is used for time-marching the solution. Firstly, the residual ($\mathbf{R e s}$) is evaluated through the inviscid and viscous spatial discretization procedure outlined in the preceding sections. Subsequently, the solution at the end of each RK sub-step is computed using the following three expressions:

\begin{equation}
    \begin{aligned}
\mathbf{U}^{(1)} &=\mathbf{U}^{\mathbf{n}}+\Delta t \mathbf{R e s}\left(\mathbf{U}^{\mathbf{n}}\right) \\
\mathbf{U}^{(2)} &=\frac{3}{4} \mathbf{U}^{\mathbf{n}}+\frac{1}{4} \mathbf{U}^{(1)}+\frac{1}{4} \Delta t \mathbf{R e s}\left(\mathbf{U}^{(1)}\right) \\
\mathbf{U}^{\mathbf{n}+\mathbf{1}} &=\frac{1}{3} \mathbf{U}^{\mathbf{n}}+\frac{2}{3} \mathbf{U}^{(2)}+\frac{2}{3} \Delta t \mathbf{R e s}\left(\mathbf{U}^{(2)}\right),
\end{aligned}
\end{equation}

\subsection{Stability criterion}
Von Neumann stability analysis was conducted using the explicit first order forward Euler time-integration \cite{hindmarsh1984stability} for one-dimensional diffusion system \cite{nishikawa2010beyond} to derive an \textit{approximate} stable time-step criterion. The resulting time-step criterion although not accurate for the present time discretization, can still be used to get rough estimates, since the time-step criterion for Euler time-integration tends to be more restrictive than that for RK-3. Furthermore, the viscous coefficient $\mu$ is considered constant in the analysis for simplicity. Following the Fourier error analysis discussed, the magnitude of amplification factor is restricted to below unity, leading to the following two conditions for ME4-Opti and ME6-Opti schemes respectively:

\begin{equation}
\begin{aligned}
    \left|1+\frac{\mu \Delta t}{\Delta x^2}\left( -\frac{558379}{180000} + \frac{142793 \cos (k)}{40000} - \frac{52793 \cos (2k)}{100000} + \frac{108379 \cos (3k)}{1800000} \right)\right| &\leq 1, \\
    \left|1+\frac{\mu \Delta t}{\Delta x^2}\left( -\frac{9223447}{2880000} + \frac{2261149 \cos (k)}{600000} - \frac{821149 \cos (2k)}{1200000} + \frac{1663447 \cos (3k)}{12600000} - \frac{461149 \cos (4k)}{33600000} \right)\right| &\leq 1.
\end{aligned}
\end{equation}

As the derived schemes demonstrate maximum damping magnitude at the grid cut-off wavenumber (Fig. \ref{optimizedSchemes}), the stability criterion can be obtained by substituting $k=\pi$ in the above two equations. This yields the following:

\begin{equation} \label{ts_restrict}
\begin{aligned}
    \text{ME4-Opti:} \quad \Delta t &\le \frac{\Delta x^2}{3.63 \mu}, \\
    \text{ME6-Opti:} \quad \Delta t &\le \frac{\Delta x^2}{3.90 \mu}
\end{aligned}
\end{equation}

The mixed derivative terms also hold approximately the same criteria. All the simulations presented in the current work are performed maintaining a Courant Friedrichs Lewy (CFL) number of 0.5, unless stated otherwise for a specific case. The time-step ($\Delta t$) satisfying the CFL condition is evaluated during each iteration as follows:

\begin{equation}
\Delta t=\mathrm{CFL} \times \min \biggr[ \underbrace{\min _{\text {cells }}\left(\frac{\Delta x}{|u|+c}, \frac{\Delta y}{|v|+c}, \frac{\Delta z}{|w|+c}\right) }_{\text{for convection terms}} , \underbrace{ \min _{\text {cells }}\left(\frac{\Delta x^2}{\mathcal{D} \mu}, \frac{\Delta y^2}{\mathcal{D} \mu}, \frac{\Delta z^2}{\mathcal{D} \mu}\right)}_{\text{for viscous terms}} \biggr],
\end{equation}

\noindent where the parameter $\mathcal{D}$ is equal to $3.63$ and $3.90$ for ME4-Opti and ME6-Opti respectively from Eqn. \ref{ts_restrict}.

%% file: results.tex
\section{Numerical tests and discussion} \label{sec:results}

\subsection{Order of accuracy - straight derivative terms} \label{sec:OA}
We test the order of accuracy of the proposed schemes by discretizing the following straight second derivative term:

\begin{equation}\label{OAstraight}
    \frac{\partial}{\partial x}\left( \mu \frac{\partial f}{\partial x} \right)
\end{equation}

This equation (Eqn. \ref{OAstraight}) precisely generalizes all the straight viscous terms that appear in the compressible Navier-Stokes equations. In particular we test the order of accuracy for a non-linear case (variable $\mu$) As described in Ref. \cite{shen2009high}, the function $f$ is specified using a high-frequency sine wave $f(x) = \sin (10x)$ over a uniform grid containing $N$ points spanning in the range $x \in [0,1]$. The viscous coefficient varies in space exponentially as $\mu(x) = 0.1 e^{2x}$. The grid is progressively refined from $20$ points to $320$ points in four steps by doubling the grid count each step. The L1 norm error is used to assess the order of accuracy.


\begin{table}[h!]
    \centering
    \begin{tabular}{lcccccccc}
        \toprule & \multicolumn{2}{c}{ Nishikawa $(2010)$ \cite{Nishikawa2010}} & \multicolumn{2}{c}{ Shen et al. $(2009)$ \cite{shen2009high} } & \multicolumn{2}{c}{ ME4-Base } & \multicolumn{2}{c}{ ME4-Opti } \\
        & L1 Error & Order & L1 Error & Order & L1 Error & Order & L1 Error & Order \\
        \midrule $\mathrm{N}=20$ & 7.62E-02 & - & 7.03E-03 & - & 3.13E-02 & - & 2.32E-02 & - \\
        $\mathrm{N}=40$ & 1.69E-02 & 2.171 & 4.36E-04 & 4.011 & 1.99E-03 & 3.978 & 1.54E-03 & 3.908 \\
        $\mathrm{N}=80$ & 4.10E-03 & 2.044 & 2.72E-05 & 4.005 & 1.25E-04 & 3.996 & 9.78E-05 & 3.981 \\
        $\mathrm{N}=160$ & 1.02E-03 & 2.011 & 1.70E-06 & 4.001 & 7.79E-06 & 3.998 & 6.14E-06 & 3.994 \\
        $\mathrm{N}=320$ & 2.54E-04 & 2.002 & 1.06E-07 & 4.000 & 4.87E-07 & 3.999 & 3.84E-07 & 3.999 \\
        \bottomrule
    \end{tabular}
    \caption{L1 errors and order of accuracy demonstrations for various 4th order straight derivative discretization schemes (the scheme of Nishikawa \cite{Nishikawa2010} is designed only for fourth order linear accuracy).}
    \label{OAStraight1}
\end{table}

\begin{table}[h!]
    \centering
    \begin{tabular}{ccccccc}
    \toprule & \multicolumn{2}{c}{ Shen \& Zha (2010) \cite{shen2010large} } & \multicolumn{2}{c}{ ME6-Base } & \multicolumn{2}{c}{ ME6-Opti } \\
    & L1 Error & Order & L1 Error & Order & L1 Error & Order \\
    \midrule 
    $\mathrm{N}=20$  & 2.73E-04 & -     & 4.02E-04 & -     & 1.51E-03 & -     \\
    $\mathrm{N}=40$  & 4.25E-06 & 6.008 & 6.37E-06 & 5.981 & 2.56E-05 & 5.886 \\
    $\mathrm{N}=80$  & 6.65E-08 & 5.999 & 1.00E-07 & 5.992 & 4.07E-07 & 5.975 \\
    $\mathrm{N}=160$ & 1.04E-09 & 6.000 & 1.57E-09 & 5.999 & 6.38E-09 & 5.995 \\
    $\mathrm{N}=320$ & 1.74E-11 & 5.897 & 2.52E-11 & 5.955 & 9.30E-11 & 6.100 \\
    \bottomrule
    \end{tabular}
    \caption{L1 errors and the corresponding order of accuracy demonstrations for various 6th order straight derivative discretization schemes.}
    \label{OAStraight2}
\end{table}

Tables \ref{OAStraight1} and \ref{OAStraight2} show the L1 norm errors and the order of accuracy obtained by different schemes. It is clear that the proposed fourth and sixth order schemes achieve the expected high order accuracy. The second order degeneracy of the Nishikawa's scheme \cite{Nishikawa2010} is due to the non-linearity of the test case. The viscous coefficient is interpolated using the standard fourth order formula (Eqn. \ref{InterPMu-base}) while testing the order of accuracy for Nishikawa's scheme. Furthermore, the absolute values of L1 errors for the ME6-Opti scheme are closely comparable to the base schemes, although slightly higher due to ME6-Opti's smaller stencil. The increase in error compared to the sixth-order scheme of Shen \& Zha \cite{shen2010large} can be attributed to the trade-off made during the optimization process between the leading truncation error, stencil size, and the spectral representation of optimized schemes in the high wavenumber range. However, the superior spectral properties of the ME6-Opti scheme are expected to provide good damping and stability characteristics in actual flow simulations which will be demonstrated subsequently in Sec. \ref{sec:DPSL}.

\subsection{Order of accuracy - mixed derivative terms} \label{sec:OA2}
To test the order of accuracy of the mixed second derivative terms, we propose to discretize the following non-linear term over a 2-D $x-y$ domain.

\begin{equation}\label{OAmixedeqn}
    \frac{\partial}{\partial x}\left( \mu \frac{\partial g}{\partial y} \right)
\end{equation}

This term generalizes all the mixed second derivative viscous terms that appear in the governing equations. Similar to the previous one dimensional case, we employ a high frequency wave and an exponentially increasing function to specify $g$ and $\mu$ respectively.

\begin{equation}
    g(x,y)= \sin [10(x+y)], \quad \mu (x,y) = 0.1 e^{2(x+y)}, \quad x \in [0,1], \quad y \in [0,1]
\end{equation}

The grid is uniformly discretized using $N$ points in both $x$ and $y$ directions.


\begin{table}[h!]
    \centering
    \begin{tabular}{lccccccccc}
    \toprule & \multicolumn{2}{c}{ Extended linear scheme} & \multicolumn{2}{c}{ Shen et al. $(2009)$ \cite{shen2009high}} & \multicolumn{2}{c}{ ME4-Base } & \multicolumn{2}{c}{ ME4-Opti } \\
     & L1 Error & Order & L1 Error & Order & L1 Error & Order & L1 Error & Order \\
    \midrule $\mathrm{N}=20$ & 5.25 & - & 2.62E-01 & - & 2.42E-01 & - & 3.47E-01 & - \\
    $\mathrm{N}=40$ & 1.36 & 1.946 & 1.68E-02 & 3.964 & 1.55E-02 & 3.968 & 2.32E-02 & 3.905 \\
    $\mathrm{N}=80$ & 0.34 & 1.986 & 1.06E-03 & 3.992 & 9.72E-04 & 3.993 & 1.48E-03 & 3.966 \\
    $\mathrm{N}=160$ & 8.62E-02 & 1.996 & 6.60E-05 & 3.998 & 6.08E-05 & 3.999 & 9.36E-05 & 3.986 \\
    $\mathrm{N}=320$ & 2.16E-02 & 1.999 & 4.13E-06 & 3.999 & 3.80E-06 & 3.999 & 5.87E-06 & 3.994 \\
    \bottomrule
    \end{tabular}
    \caption{L1 errors and the corresponding order of accuracy demonstrations for various 4th order mixed derivative discretization schemes.}
    \label{OAmixed1}
\end{table}

\begin{table}[h!]
    \centering
    \begin{tabular}{ccccccc}
    \toprule & \multicolumn{2}{c}{ Shen \& Zha (2010) \cite{shen2010large}} & \multicolumn{2}{c}{ ME6-Base } & \multicolumn{2}{c}{ ME6-Opti } \\
    & L1 Error & Order & L1 Error & Order & L1 Error & Order \\
    \midrule 
    $\mathrm{N}=20$  & 1.36E-02 & -     & 1.21E-02 & -     & 2.43E-02 & - \\
    $\mathrm{N}=40$  & 2.21E-04 & 5.945 & 1.96E-04 & 5.949 & 4.13E-04 & 5.877 \\
    $\mathrm{N}=80$  & 3.49E-06 & 5.987 & 3.10E-06 & 5.988 & 6.65E-06 & 5.959 \\
    $\mathrm{N}=160$ & 5.46E-08 & 5.997 & 4.84E-08 & 5.997 & 1.05E-07 & 5.985 \\
    $\mathrm{N}=320$ & 8.54E-10 & 5.998 & 7.58E-10 & 5.997 & 1.65E-09 & 5.993 \\
    \bottomrule
    \end{tabular}
    \caption{L1 errors and the corresponding order of accuracy demonstrations for various 6th order mixed derivative discretization schemes.}
    \label{OAmixed2}
\end{table}

Tables \ref{OAmixed1} and \ref{OAmixed2} show the L1 norm errors for various schemes for the test function $g(x,y)$. All the tested high order schemes can be noted to achieve their respective desired order of accuracy. We also test the fourth order linearly accurate extended version of Nishikawa's scheme for mixed derivatives as employed in Ref. \cite{chamarthi2022importance} through the $\alpha$-damping approach. The corresponding results are presented in Table \ref{OAmixed1} under the name `Extended linear scheme'. The accuracy for this linear scheme was noted to be comparatively low.\\
The ME6-Base scheme demonstrates the minimum absolute error magnitude amongst all the schemes compared. The Optimized schemes show a slightly higher error than the schemes of Shen et al. \cite{shen2009high} and Shen \& Zha \cite{shen2010large}, due to the extra filtering term added in to Eqn. \ref{eqn:InterpFilter}, which is the trade-off made in exchange for high wavenumber accuracy. Despite the filtering term, the errors associated with the optimized schemes are still close and comparable to that of other schemes at all grid resolutions.

\subsection{Non-linear diffusion equation} \label{sec:Diffusion}
The impact of non-linear high-order accuracy is assessed in the present example using the unsteady non-linear diffusion equation. The Eqn. \ref{non-lin-diff} under consideration includes the straight non-linear diffusion term, a feature shared with the compressible Navier-Stokes and also non-linear Burgers equations.

\begin{equation}\label{non-lin-diff}
    \frac{\partial f}{\partial t} = \frac{\partial}{\partial x}\left(\nu \frac{\partial f}{\partial x}\right)
\end{equation}

The problem is solved in a periodic domain spanning the range $x \in [0, 1]$. The field variable $f$ is initialized with a high-frequency wave as $f(x,t=0)= \sin (16 \pi x)$. The diffusivity coefficient $\nu$ is fixed to be equal to $\nu (x)= \cos (16 \pi x)$. The problem is solved until an end time of $t_{end}=0.0025$ on a uniform grid containing 144 points accounting 8 cycles per unit length and an effective grid wavenumber of $k_{grid}=\pi/4.5$. A constant global time-step of $\Delta t=1.2 \times 10^{-6}$ is used. The reference solution is computed on a grid containing 1024 points ($k_{grid}=\pi/32$) using the ME4-Base scheme, at which grid count the solution is found to have sufficiently converged for the purpose of present comparison. \\

The high-order accurate schemes developed in this study are compared to the fourth-order linearly accurate scheme proposed by Nishikawa \cite{nishikawa2010beyond} (designed for handling constant viscosity cases). It's crucial to highlight that this seemingly disproportionate comparison is conducted solely to assess the magnitude of error introduced into the solution when the scheme lacks non-linear high-order accuracy. The discretization formulation for Nishikawa's scheme \cite{nishikawa2010beyond} employed here is briefly detailed in Sec. \ref{sec:intro}, where $\alpha = 8/3$ is chosen to achieve fourth-order accuracy. The viscosity at the interface is computed by averaging viscosity values from the two nodes situated on either side of each interface. 


\begin{figure}[h!]
    \centering
    \includegraphics[width=\textwidth]{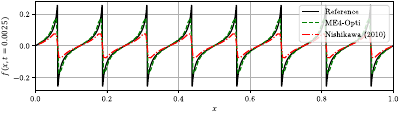}
    \caption{Profile of `$f$' evaluated using the ME4-Opti and the fourth order linear scheme of Nishikawa \cite{nishikawa2010beyond} at $t=0.0025$ shown in comparison with the reference solution.}
    \label{fig:non-linDiff1}
\end{figure}

\begin{figure}[h!]
    \centering
    \includegraphics[width=\textwidth]{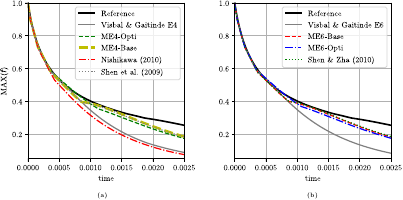}
    \caption{Temporal evolution of the solution's global maximum computed using various (a) 4th order and (b) 6th order schemes.}
    \label{fig:non-linDiff2}
\end{figure}

Fig. \ref{fig:non-linDiff1} presents the solution at $t=0.0025$. The peak value for the solutions corresponding to reference, ME4-Opti, and Nishikawa's linear scheme \cite{nishikawa2010beyond} can be noted to be equal to $0.257$, $0.174$, and $0.075$ respectively. The linear scheme of Nishikawa was noted to over dissipate the solution in this scenario, and the solution of ME4-Opti was closer to the reference value. \\

The evolution of the peak value, MAX($f$), over time is depicted in Fig. \ref{fig:non-linDiff2} for various fourth and sixth-order schemes. While very close to the solutions of ME4-Opti and Shen et al. \cite{shen2009high}, the most accurate among the compared fourth-order schemes was the solution of the ME4-Base scheme \cite{zingg2000comparison, de1999aerodynamic}. The slight difference between ME4-Opti and ME4-Base can be attributed to the truncation error penalty resulting from the optimization of spectral properties. A similar trend was observed for the sixth-order schemes. On the other hand, despite the Visbal E4 scheme being high-order, the error corresponding to its solution is observed to be relatively high. Additionally, spurious oscillations were also observed in its solution field. This is speculated to be attributed to it's odd-even decoupling nature.

\subsection{Doubly Periodic Shear Layer} \label{sec:DPSL}
The case of the Doubly Periodic Shear Layer \cite{bell1989second, clausen2013entropically,achu2021entropically} consists of a stream of fluid moving through a stationary field of uniform pressure and density. The flow evolves to form two shear layers on the top and bottom interfaces. The flow evolves to generate two shear layer vortices of opposing spin on the top and bottom. The Reynolds number and Mach number of the flow are set to $10,000$ and $0.1$, respectively. The simulation is conducted in a periodic domain of size $[0,1]^2$ until an end time of $t=1$.\\

\textcolor{magenta}{The flowfield comprises large amplitude spatial gradients in velocity across the shear layer, where viscous diffusion, grid size, and the discretization of advection terms play crucial roles in the formation of vortex structures, as demonstrated in the works of Brown \cite{brown1995performance}, Minion \& Brown \cite{minion1997performance}, and Drikakis \& Smolarkiewicz \cite{drikakis2001spurious}. In this context, it is logical to mask the effect of inviscid scheme and isolate the effect of various viscous flux discretization schemes. Therefore, a fixed inviscid flux discretization scheme has been employed for all cases presented in the present section.} The inviscid fluxes at each node are first computed, and then a standard explicit sixth-order central scheme \cite{lele1992compact, visbal1999high} is employed to compute the derivatives of flux field directly. Although the dissipation error is zero (since central), the dispersion error of the scheme can lead to an accumulation of energy, resulting in unphysical oscillations in the flowfield that grow over time and, consequently, numerical instability. While the natural viscosity of the Navier-Stokes equations helps dampen these oscillations, sufficient damping is required, especially in the high wavenumber region, to maintain stability at high Reynolds numbers, such as in the present case. To address this, the conserved variables are spatially filtered \cite{Visbal2002, vadlamani2018distributed, delorme2017simple} using the optimized sixth-order explicit filter of Bogey et al. \cite{bogey2009shock} after time integration after every fixed number of time-integration cycles. The filtering operation relaxes energy from scales close to the grid cutoff wavenumber, counteracting the effect of the dispersion error from the inviscid scheme and stabilizing the solution. The filtering cycle, denoted as $\Theta$, is defined as the number of time-steps per which one filtering operation is performed, and it is being varied to assess the simulation stability of various schemes listed in Table \ref{ResolvingEff} due to the differences in their spectral resolution.\\

\begin{figure}[h!]
    \centering
    \includegraphics[width=0.90\textwidth]{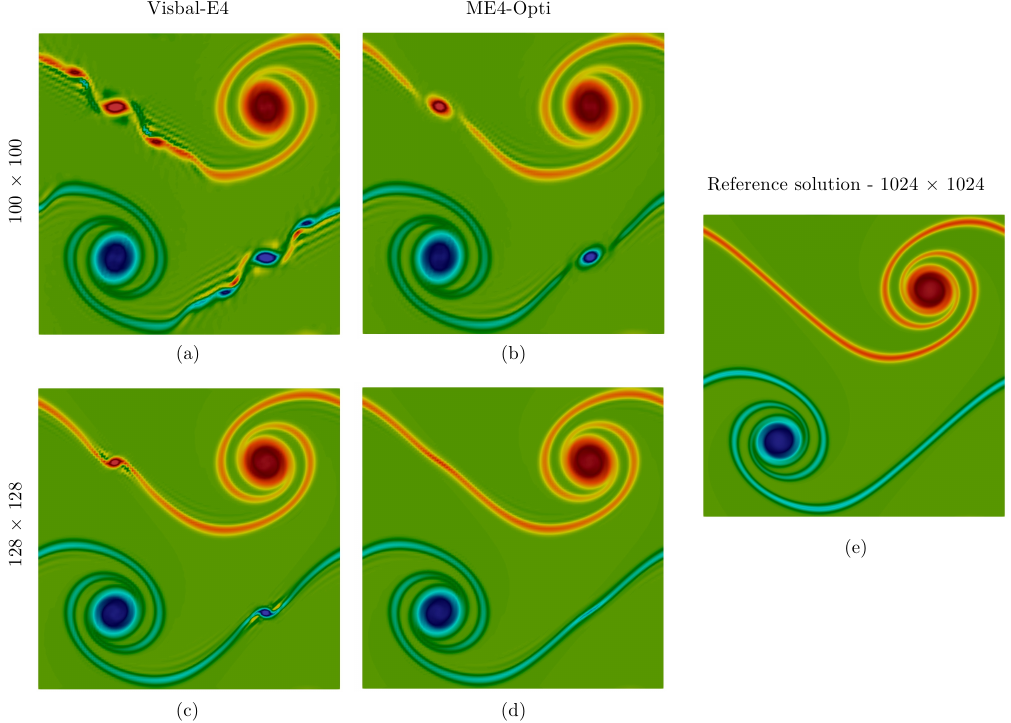}
    \caption{(a, b, c, d) Comparison of z-vorticity fields simulated using the Visbal \& Gaitonde E4 and ME4-Opti viscous schemes at grid resolutions of $100 \times 100$ and $128 \times 128$. The filtering cycle ($\Theta$) was maintained at one filter per 140 time-steps and 190 time-steps for the $100 \times 100$ and $128 \times 128$ grids, respectively. The braid vortex is circled in (a). The reference solution from a fine grid simulation using the Visbal E4 scheme is shown on the side in (e).}
    \label{fig:DPSL1}
\end{figure}

Fig. \ref{fig:DPSL1} compares the $x$-vorticity fields ($\Omega_x$) of the solution at $t=1$, simulated using the fourth-order Visbal \& Gaitonde E4 scheme and the present ME4-Opti scheme. The test case produces a pair of small braid vortices on the top and bottom shear layers when the grid resolution and/or dissipation are not sufficient enough. The size of the braid vortex (circled in Fig. \ref{fig:DPSL1}a) can vary depending on the amount of dissipation being added, both from the viscous scheme and the filtering operation. Since the filtering cycle $\Theta$ is kept constant, the dissipation from filtering remains fixed for both simulations, and the difference in the solution is solely attributed to the choice of the viscous scheme.\\

The vorticity fields in Fig. \ref{fig:DPSL1} show that the braid vortex size of the solution corresponding to Visbal \& Gaitonde E4 scheme is notably larger and significantly off from the reference solution. This observation holds true for both $100 \times 100$ and $128 \times 128$ grids. This emphasizes the influence of good high wavenumber spectral resolution of the proposed schemes. It is worth noting that, although the braid vortex would vanish at a low enough value of $\Theta$ for both schemes, the presented filtering cycle ($\Theta$) values, which are $140$ and $190$ for $100 \times 100$ and $128 \times 128$ grids, respectively, have been deliberately chosen to highlight the differences that the spectral properties of the schemes show on the solution. This demonstrates that with a fixed inviscid flux discretization scheme, the viscous scheme with better spectral properties converges more rapidly to the actual solution with grid refinement.\\

\begin{figure}[h!]
    \centering
    \includegraphics[width=\textwidth]{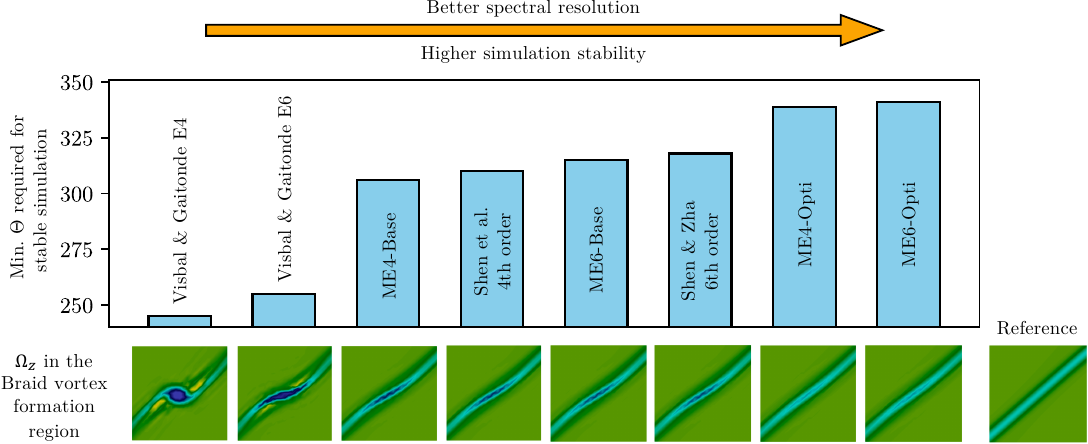}
    \caption{Variation in the minimum filter cycle ($\Theta$) value required to maintain a stable simulation for different schemes is observed for the doubly periodic shear layer test case, at a grid resolution of $128 \times 128$. Any value of $\Theta$ lower than the value displayed by the bar plot will result in an unstable solution and numerical blow-up. The decay of the braid vortex concerning the improvement in the spectral resolution of the scheme is depicted at the bottom.}
    \label{fig:DPSL2}
\end{figure}

To assess the simulation stability against the aliasing errors from the inviscid central discretization scheme (not to be confused with the time-step restriction) versus the spectral property of the viscous flux discretization scheme, we present the minimum number of iterations per which filtering (minimum $\Theta$) is to be performed to maintain a stable simulation. A relatively more stable scheme would require lesser amount of filtering and thus high $\Theta$. Fig. \ref{fig:DPSL2} shows the variation of $\Theta$ for various viscous flux discretization schemes obtained on a $128 \times 128$ grid. The size of the braid vortex for the corresponding simulation with a fixed value of $\Theta=140$ is shown on the bottom. The braid vortex can be seen to vanish as the resolution of the scheme is improved. This is attributed to the accurate viscous damping provided by the scheme. The results are in consistent with the resolving efficiency ($e_v$) and spectral viscosity ($\nu_s^{''}$) values reported in Table \ref{ResolvingEff}. 

\subsection{Kelvin-Helmholtz instability} \label{KHI}
Next, we simulate the Kelvin-Helmholtz instability \cite{ryu2000magnetohydrodynamic, san2015evaluation}, a prevalent flow phenomenon observed in various fluid flows featuring shear layers \textcolor{magenta}{and vortical structures \cite{tsoutsanis2015comparison}}. The flow is initialized as follows:

\begin{equation} \label{KHini}
    \begin{aligned}
    &p=2.5, \quad \rho(x, y)=\left\{\begin{array}{l}
    2, \text { If } 0.25 \leq y \leq 0.75 \\
    1, \text { else }
    \end{array}\right. \\
    &u(x, y)=\left\{\begin{array}{l}
    0.5, \quad \text { If } 0.25 \leq y \leq 0.75 \\
    -0.5, \text { else },
    \end{array}\right. \\
    &v(x, y)=0.1 \sin (4 \pi x)\left\{\exp \left[-\frac{(y-0.75)^{2}}{2 \sigma^{2}}\right]+\exp \left[-\frac{(y-0.25)^{2}}{2 \sigma^{2}}\right]\right\}, \text { where } \sigma=0.05/ \sqrt{2}.\
    \end{aligned}
\end{equation}

The flow is simulated within a two-dimensional periodic domain of size $[0,1]\times[0,1]$. The shear layer is perturbed to initiate the instability using the aforementioned mentioned $y$-velocity expression. Simulations were conducted at a Reynolds number of Re=$200$ and a Mach number of Mach=$0.1$, with a grid resolution of $1100 \times 1100$. Similar to the Doubly Periodic Shear Layer case discussed in Sec. \ref{sec:DPSL}, the inviscid fluxes in the current test case were discretized using the standard sixth-order explicit central scheme. Conservative variables were filtered through the sixth-order low-pass filter of Bogey et al. \cite{bogey2009shock, bogey2011finite} after every 450 time iterations. This choice of $\Theta=450$ was implemented to maintain simulation stability without introducing dissipation from the filtering process, thus isolating the effect of viscous schemes. A constant non-dimensional time-step of $\Delta t = 1.5 \times 10^{-4}$ was employed, and the simulation ran until reaching an end-time of $t=1$.

\begin{figure}[h!]
    \centering
    \includegraphics[width=\textwidth]{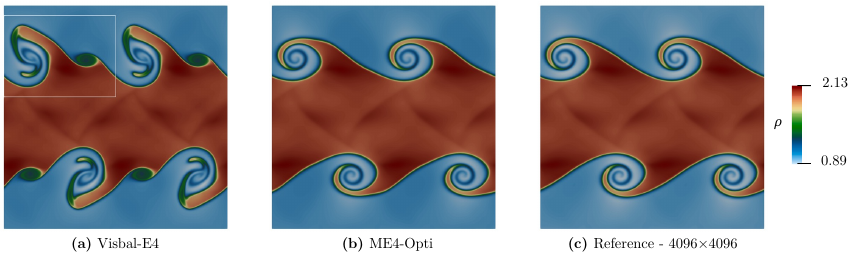}
    \caption{Density contours for the fourth order Visbal-E4 and ME4-Opti schemes in comparison with the reference solution.}
    \label{fig:KHI1}
\end{figure}

\begin{figure}[h!]
    \centering
    \includegraphics[width=\textwidth]{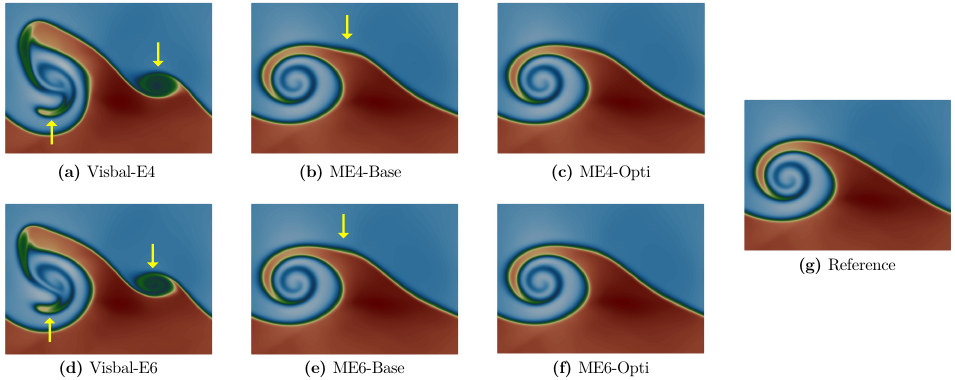}
    \caption{Density contours near the vortex structure obtained from various schemes. The presented contours are presented from a cropped-in portion of the computational domain represented by the white box in Fig. \ref{fig:KHI1}(a).}
    \label{fig:KHI2}
\end{figure}

Figure \ref{fig:KHI1} displays the density field of the Kelvin-Helmholtz instability obtained using the fourth-order Visbal-E4 scheme \cite{visbal1999high}, the current ME4-Opti, and the converged reference solution obtained on a $4096 \times 4096$ grid using the Visbal-E4 scheme. The solution qualitatively demonstrates the accuracy of the present ME4-Opti scheme in comparison to the reference solution. Conversely, the Visbal-E4 scheme appears to yield an under-dissipated solution with flow instabilities that are absent in the reference solution. Furthermore, Fig. \ref{fig:KHI2} compares the solutions obtained using other schemes, including the sixth-order ones. The arrows in Fig. \ref{fig:KHI2} denote the locations of the onset of flow instability. It can be inferred that an increase in the spectral accuracy of the scheme leads to the appropriate amount of viscous dissipation, resulting in a more accurate solution.

\begin{figure}[h!]
    \centering
    \includegraphics[width=\textwidth]{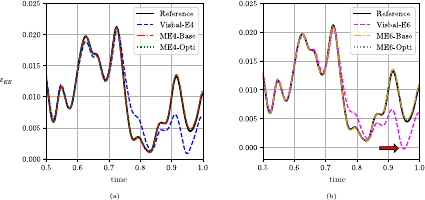}
    \caption{Kinetic energy dissipation rate $\varepsilon_{KE}$ histories between $t=0.5$ and $1$ for various (a) fourth and (b) sixth order schemes.}
    \label{fig:KHI3}
\end{figure}

To also present a quantitative comparison, the kinetic energy dissipation rate in the flow is tracked in time. The normalized cumulative kinetic energy of the domain is computed using Eqn. \ref{KE-cum}. The kinetic energy dissipation rate ($\varepsilon_{KE}$) is computed by taking a derivative of KE (Eqn. \ref{KE-der}). In Fig. \ref{fig:KHI3}, the evolution of the kinetic energy dissipation rate is shown, obtained using both fourth and sixth order schemes. While the present base and optimized Midpoint-based Explicit schemes (ME4 and ME6) have converged to the reference solution, there is a noticeable discrepancy for the low-resolution Visbal-E4 and Visbal-E6 schemes. Additionally, an undershoot below the value of $\varepsilon_{KE}=0$ is observed in the plot of the Visbal-E6 scheme at $t\approx0.95$, indicating an unphysical increase in kinetic energy within the system.

\begin{equation} \label{KE-cum}
    KE = \frac{1}{\Delta x \Delta y} \sum_{cells} \frac{1}{2} \rho (u^2+v^2)
\end{equation}

\begin{equation} \label{KE-der}
     \varepsilon_{KE} = -\frac{\partial (\text{KE})}{\partial t}
\end{equation}

\begingroup
\color{blue}
\subsection{Odd-even decoupling} \label{sec:odd-even}
The Odd-Even decoupling test first proposed by Quirk in \cite{quirk1997contribution} was introduced to evaluate the robustness of numerical schemes in handling grid-aligned shocks and discontinuities without developing spurious oscillations. The test case was originally proposed to test inviscid schemes for their ability to prevent odd-even decoupling \cite{gallice2022entropy}. Inviscid schemes that fail to produce a clean shock in this test will exhibit saw-tooth type oscillations behind the shock wave. However, in this section, we extend this test case to evaluate the odd-even decoupling-free property of the proposed viscous flux discretization schemes, ME4-Opti and ME6-Opti.


The test case involves a shock wave of Mach 6 propagating in a computational domain of size $x \in [0,800]$ and $y \in [0,20]$ from left to right. The shock wave is initialized at the center of the domain at $x = 400$. The grid is uniformly discretized with unit spacing making $\Delta x = \Delta y = 1$ or the number of cells $N_x= 800$ and $N_y = 20$. However, the grid points along the horizontal center line, i.e., along $y = 10$ are slightly perturbed as $y \leftarrow y + (-1)^{i}10^{-6} \quad \forall, i \in [0,N_x]$, to promote odd-even decoupling phenomenon along the length of the shock. The right half of the domain was initialized with quiescent pre-shock flow conditions: $\left[\rho, u, v, p\right]_{\text{R}} = \left[\gamma, 0, 0, 1\right]$. The post-shock conditions on the left half of the domain can be computed using the Rankine-Hugoniot relations:

\begin{equation}
    \left[\rho, u, v, p\right]_{\text{L}} = \left[
        \frac{\gamma (\gamma+1) \mathrm{M}^2}{(\gamma-1) \mathrm{M}^2+2}, \,
        \frac{2\left(\mathrm{M}^2-1\right)}{(\gamma+1) \mathrm{M}}, \,
        0, \,
        \frac{2 \gamma \mathrm{M}^2-(\gamma-1)}{(\gamma+1)}
    \right],
\end{equation}

\noindent where M=6. The simulation was run until an end time of $t_{\text{end}} = 50$. Given that Quirk's test case \cite{quirk1997contribution} was originally designed to evaluate inviscid schemes, the following modifications were made to evaluate the present viscous flux discretization schemes:

\begin{itemize}
    \item The inviscid setting of the originally proposed test case \cite{quirk1997contribution} was modified to a viscous setting, including all the viscous terms described in the governing equations section (Sec. \ref{sec:gov-eqns}). A Reynolds number of Re = 1000 was considered.
    \item The top and bottom boundary conditions of the original Quirk's test case \cite{quirk1997contribution} were changed to periodic instead of walls. This prevents the formation of boundary layers, which could potentially lead to shock wave-boundary layer interaction phenomena and interfere with our observations.
    \item For the inviscid flux discretization, we used a standard WENO-5 scheme \cite{Jiang1995}, alongside the HLL Riemann solver \cite{toro2009riemann} and primitive variable reconstruction \cite{kakumani2022use}. We chose WENO-5 with HLL because it inherently avoids odd-even decoupling that may arise from inviscid flux discretization \cite{kumar2013weno}. Therefore, any oscillations observed behind the shock at \( t_{\text{end}} = 50 \) will be purely due to the viscous flux discretization.
\end{itemize}

\begin{figure}[h!]
    \centering
    \includegraphics[width=\textwidth]{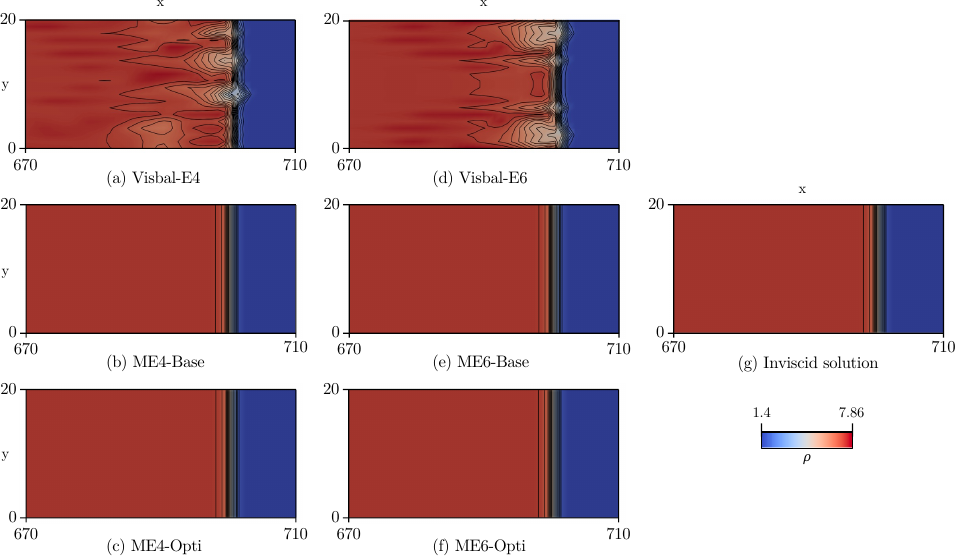}
    \caption{Density field at \( t = 50 \) obtained using various fourth and sixth order viscous flux discretization schemes (a-f) and the reference solution obtained using an inviscid simulation (g). The density iso-lines are plotted using twenty uniformly spaced values in the range [1.4, 7.86].}
    \label{fig:oddEven}
\end{figure}

The test was conducted using two classes of viscous schemes. The first class includes schemes that are prone to odd-even decoupling phenomenon, specifically Visbal-E4 and Visbal-E6. The second class comprises schemes that are free from odd-even decoupling, such as the ones proposed in this paper, ME4-Base and ME4-Opti. Other schemes that are also free from odd-even decoupling, as discussed in section \ref{sec:intro}, include those by Shen et al. \cite{shen2009high, shen2010large}, Zingg et al. \cite{zingg2000comparison} and De Rango \cite{de1999aerodynamic}.

Figure \ref{fig:oddEven} shows the results obtained at \( t = 50 \) using both classes of schemes. It can be seen that the Visbal E4/E6 results exhibit oscillations behind the shock, while the ME4-Base/Opti and ME6-Base/Opti show clean results similar to that of inviscid case shown in Fig. \ref{fig:oddEven}(g). Although all the schemes tested here are of high-order accuracy, the oscillations observed in the first class of schemes are attributed to their odd-even decoupling nature. This demonstrates that preventing odd-even decoupling through appropriate choice of viscous flux discretization scheme is crucial for avoiding oscillations in this test case, and these conclusions can be extrapolated to other flows with strong shocks. 
\endgroup

\subsection{Mach 1.56 impinging choked supersonic jet} \label{sec:SupersonicJet}
The primary purpose of this test case is to demonstrate the working of the presented optimized viscous flux discretization schemes on non-uniform curvilinear grids.  Extending the present ME4-Opti and ME6-Opti to curvilinear geometries is straightforward, and the specifics concerning metric computation and transformed fluxes can be adapted from those outlined extensively in Ref. \cite{chandravamsi2023application}. For the sake of brevity they are not repeated here. In the present sixth order accurate simulation, the viscous fluxes are discretized using the ME6-Opti scheme. The inviscid fluxes on the other hand are discretized using a sixth order eleven point stencil Monotonocity Preserving Optimized Upwind Reconstruction Scheme (MP-OURS6), a modification to the standard Monotonocity Preserving (MP) scheme of Suresh and Huynh \cite{suresh1997accurate}. The details regarding the inviscid scheme can be found in \ref{sec:OURS-6}. Resolving the supersonic jet noise \cite{edgington2019aeroacoustic} is a particularly challenging task and depends on various important factors such as: Numerical scheme, grid resolution, nozzle inlet conditions, boundary layer profile, and near wall turbulence. It is crucial to capture various flow features such as shear-layers, acoustic waves, shocks, and expansion fans with sufficient fidelity \cite{gojon2019antisymmetric, zhang2002broadband, kakumani2023gpu, zhao2000effects} to realistically resolve the noise spectrum.

\begin{figure}[h!]
    \centering
    \includegraphics[width=\textwidth*3/4]{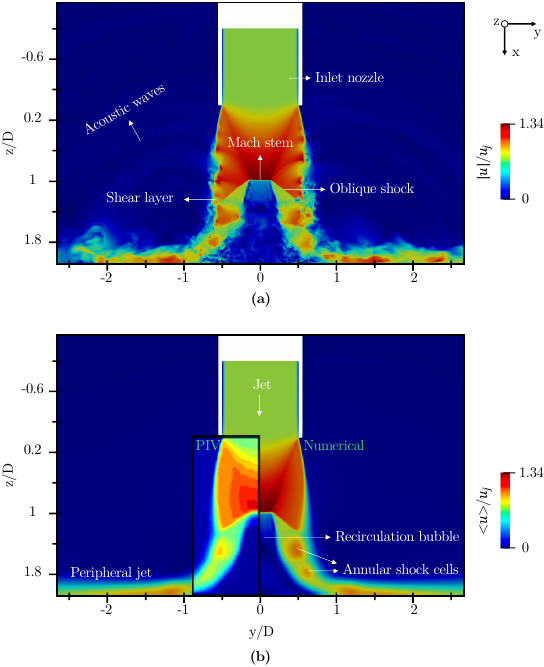}
    \caption{(a) Instantaneous velocity magnitude field of the 3-D jet impinging on a flatplate. (b) Mean Velocity magnitude field of the jet compared against the PIV data of Henderson et al. \cite{henderson2005experimental} shown by the black box. The flow simulations contain $\approx 120$ million grid points.}
    \label{fig:domains1}
\end{figure}

\subsubsection{Case setup}
The computational setup involves a 3-D choked round jet operating under a Nozzle Pressure Ratio $4.03$ and a Temperature Ratio of unity. The flow exits through a straight circular duct of diameter $D$ and impinges on a flat plate situated $2.08D$ away from the nozzle exit. The lip thickness of the exit duct is equal to $0.05D$. The ambient atmosphere is maintained at a temperature and pressure equal to $T_{\infty}=293$K and $p_{\infty}=101325$Pa. The grid and the computational domain are cylindrical and span up to a distance of $-15D$ upstream of nozzle exit and $25D$ radially outward. The ideally expanded jet Mach number and jet velocity of the flow are equal to $1.56$ and $440$ m/s respectively. The jet velocity based Reynolds number is maintained equal to Re=$\frac{D u_j}{\nu_j}=60000$. The results of this test case are being compared to the Particle Image Velocimetry (PIV) data of Henderson et al. \cite{henderson2005experimental} obtained for the same geometric configuration and at a higher Reynolds number. 

The domain was discretized uniformly using $512$ cells along the azimuthal direction. Along the radial direction, the minimum grid spacing was maintained at the nozzle lip-line with a value of $0.00375D$. The radial grid spacing within the range of $2.5D \geq r \geq 5.0D$ was kept at $\Delta r = 0.03D$ (which ideally captures waves of up to a Strouhal number $St = fD_j/u_j = 5.3$), similar to the grid utilized in the study by Gojon and Bogey in \cite{gojon2017flow}. In the streamwise direction, a minimum spacing of $\Delta x = 0.00375D$ was maintained at the nozzle exit and the impingement wall. The grid spacing remained at $\Delta x = 0.075D$ within the range $0.9D \geq x \geq 1.6D$. At the inlet of the duct at $x=-D$, a Blasius laminar boundary layer profile was imposed with a thickness of $\delta = 0.05D$. The Crocco-Busemann relation was used to specify the density profile at the inlet. No annular fluctuations were added inside the nozzle. Sponge zones are used at the farfield boundaries to avoid unwanted reflections from entering inside the computational domain. 

\subsubsection{Mean flowfield}

\begin{figure}[h!]
    \centering
    \includegraphics[width=\textwidth*3/4]{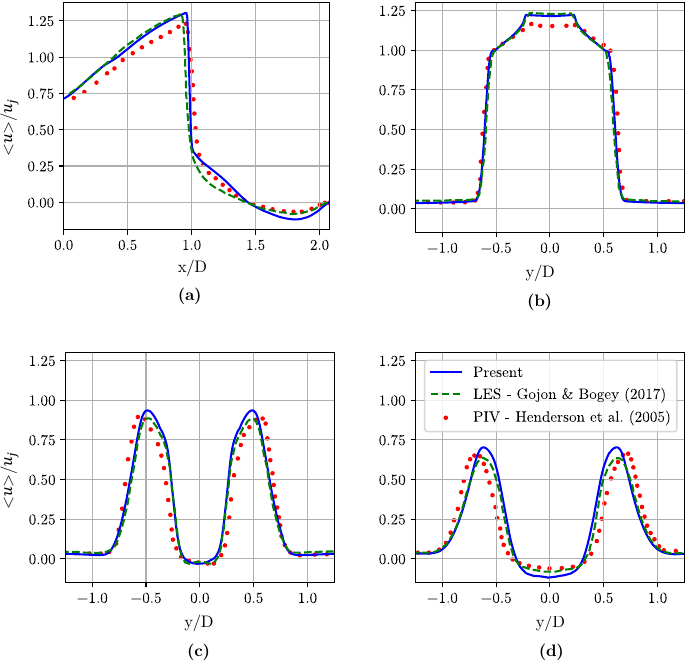}
    \caption{Comparison of present velocity data against the PIV results of Henderson et al. and LES of Gojon and Bogey \cite{gojon2017flow} along (a) Jet axis, (b) radial direction at $x/D=0.75$, (c) radial direction at $x/D=1.5$, and (d) radial direction at $x/D=1.8$.}
    \label{fig:domains2}
\end{figure}

The flow was observed to reach a statistically steady state at $t a_{\infty}/D=500$, where $a_{\infty}=343$ m/s is the ambient speed of sound. The flow statistics were then collected between the time interval $t a_{\infty}/D=500$ to $1000$. Fig. \ref{fig:domains1} presents the instantaneous and time averaged velocity magnitude flowfields of the jet. Various flow features such as shock waves, expansion fans, shear layer, recirculation bubble, and fine features such as annular shocks can be noted to have resolved clearly. A good qualitative match of the PIV data of Henderson et al. \cite{henderson2005experimental} and the present simulation is noted.\\

Fig. \ref{fig:domains2} shows the time averaged $x$-velocity profiles probed along various axes. Fig. \ref{fig:domains2}(a) features the data collected along the jet axis between the nozzle exit and impinging plate. Fig. \ref{fig:domains2}(b,c,d) shows the data collected along various radial locations along the $y$-axis direction in the domain between the nozzle exit and the impinging plate. This quantitative comparison between the experimental measurements, previous Large Eddy Simulations of Gojon and Bogey \cite{gojon2017flow} and the present simulations clearly shows a good correlation between them entrusting the present viscous flux discretization to be applicable for such practical problems as well on non-uniform meshes.



%% file: conclusions.tex
\section{Conclusions and future work}

A set of odd-even decoupling free high-resolution fourth and sixth-order viscous flux discretization schemes have been proposed to discretize both straight and mixed second derivative terms present in the Compressible Navier-Stokes equations. The schemes discussed in this paper employ velocity and temperature gradients computed at the \textit{midpoint locations} ($i+\frac{1}{2}$) directly, to eliminate the odd-even decoupling effect and result in a damping effect in the high wave number region.

Firstly, two baseline schemes, ME4-Base (fourth order), and ME6-Base (sixth order), were presented based on the previously available literature \cite{gaitonde1998high,de2021high}. These baseline schemes utilize a centered sliding stencil to compute the midpoint derivatives. The baseline schemes were then optimized for improved spectral resolution in the high wavenumber region without changing the order of accuracy or increasing the stencil size. The principle behind achieving better spectral resolution is to use the entire stencil's information to compute each mid-point derivative involved in the high-order polynomial of the second derivative (Eqn \ref{gen2ndDer}). This is unlike the sliding stencil employed in the baseline schemes (ME4-Base and ME6-Base), where only partial information from the scheme's stencil is used to compute each midpoint derivative. The proposed schemes were analyzed in spectral space, and the optimized schemes were noted to follow the exact solution closer than the base and other schemes from the literature \cite{visbal1999high,shen2010large}, resulting in higher resolving efficiency and minimum spectral viscosity. The optimization procedure has also allowed the implementation of a smaller stencil (two points smaller) for the ME6-Opti scheme which is shown to exhibit improved spectral properties compared to the Base scheme, albeit with a marginal increase in the truncation error.

We also obtain a high wavenumber damping feature in estimating the mixed derivative terms by introducing a filter term into the midpoint derivative calculation, which is absent in the baseline and other previously proposed schemes \cite{visbal1999high,shen2009high}. A summary of various properties of the second derivative discretization and the satisfaction of these properties by various schemes is presented in the following Table \ref{props}.

\begin{table}[h!]
    \centering
    \caption{Summary of various second derivative properties and their satisfaction by different schemes.}
    \includegraphics[width=\textwidth]{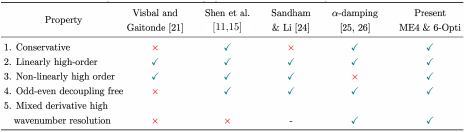}
    \label{props}
\end{table}

Numerical tests were performed to assess the order of accuracy and resolving capability of the proposed schemes. Order of accuracy tests, conducted on both straight and mixed derivative terms, demonstrate the expected fourth and sixth-order accuracy. The results from the Double Periodic Shear Layer and Kelvin-Helmholtz instability test case reveal that the proposed optimized schemes, which possess better spectral properties, are relatively more stable and produce an accurate solution that closely matches the reference solution obtained from a fine grid simulation; thus potentially reducing grid resolution requirements. The Mach 1.56 choked impinging supersonic jet, simulated at Re = 60000, was used to test the applicability of the proposed optimized schemes on non-uniform curvilinear domains. The computed mean flowfield was noted to be in good agreement with the experimental measurements and previous large eddy simulations from the literature. 

While the present paper primarily focuses on solving the Navier-Stokes equations, the second and mixed derivative schemes presented herein can also find utility in discretizing other classes of partial differential equations that involve second derivative terms, such as the Viscous Burgers equations \cite{burns1998numerical,mittal2011numerical}, Hamilton-Jacobi equation \cite{tucker2003differential,kakumani2022use} and the Poisson equation. It is worth noting that the methods proposed and tested in this paper are specific for single-component flows. However, for multi-component flows, where steep diffusion coefficient gradients are present, the robustness of the present schemes is not yet tested. The extension of the present scheme to such flows is being presented in a follow-up work. An interactive Python widget to control the spectral properties of the proposed schemes are attached to the present manuscript as supplementary material.

%% file: acknowledgements.tex
\section*{Supplementary material} \label{sec:supply}
\begin{enumerate}
    \item \textbf{Interactive spectral properties adjustment widget:}: \textcolor{blue}{\href{https://github.com/hemanthgrylls/supplementaryMaterialJCP.git}{Github repository}}
\end{enumerate}

\section*{Acknowledgement}
HC gratefully acknowledges the financial support from Technion - Israel Institute of Technology during the course of present work. The authors would like to thank the anonymous reviewer for suggesting the odd-even decoupling test case presented in Sec. \ref{sec:odd-even}.

%% file: appendix.tex
\appendix
\begingroup
\color{magenta}
\section{Equivalent numerical dissipation and Reynolds number of discretization error} \label{app:eqivalRe}
\subsection{Equivalent numerical viscosity}
The truncation errors from the second derivative schemes discussed in this work result in artificial dissipation, which is addressed in this appendix. Specifically, we examine the discretization of the diffusion term, represented as $\nu \frac{\partial }{\partial x}\left(\frac{\partial \phi}{\partial x} \right)$. Discretization can lead to either negative or positive dissipation effects. To quantify these errors, spectral viscosity \cite{lamballais2011straightforward}, $\nu^{\prime \prime}_s$, related to discretization error is introduced into the diffusion term. Spectral viscosity is defined as the amount of correction viscosity (can be either positive or negative) to be added to the actual viscosity $\nu$ to account for the discretization error. Considering the exact and numerically modeled approximation with $\phi = e^{ikx}$, the expressions after simplification will lead to the following:


\begin{equation}
    \text{Exact model:} \quad \nu \left( \frac{\partial^2 \phi}{\partial x^2} \right)_{\text{exact}} \hspace{0.2cm} = \nu \phi (-k^2),
\end{equation} \label{spec_visc1}

\begin{equation}
    \text{Numerical model:} \quad \nu \left( \frac{\partial^2 \phi}{\partial x^2} \right)_{\text{numerical}} = \nu \phi (k_{xx}^*).
\end{equation} \label{spec_visc}

\noindent Correcting the numerical model by incorporating the spectral viscosity, the diffusion term can be equivalently re-written as:

\begin{equation}
    \begin{aligned}
       \nu \left( \frac{\partial^2 \phi}{\partial x^2} \right)_{\text{numerical}} & = (\nu + \nu^{\prime \prime}_s) \left( \frac{\partial^2 \phi}{\partial x^2} \right)_{\text{exact}}, \\
        \text{i.e.,} \quad \nu \phi (k_{xx}^*) & = (\nu + \nu^{\prime \prime}_s) \phi (-k^2).
    \end{aligned}
    \label{spec_visc2}
\end{equation}

\noindent Rearranging the above terms will give the following expression for the spectral viscosity:

\begin{equation}
    \nu^{\prime \prime}_s = -\nu \frac{k_{xx}^* + k^2}{k^2}
\end{equation}

The spectral viscosity parameter, $\nu^{\prime \prime}_s$, helps understand the artificial dissipation characteristics in the wavenumber domain corresponding to various second derivative discretization schemes. Utilizing the modified wavenumber analysis from section \ref{sec:Viscmethod}, spectral viscosity profiles for various schemes are plotted as shown in Fig. \ref{spectral_visc}.

\begin{figure}[h!]
    \centering
    \includegraphics[width=150mm]{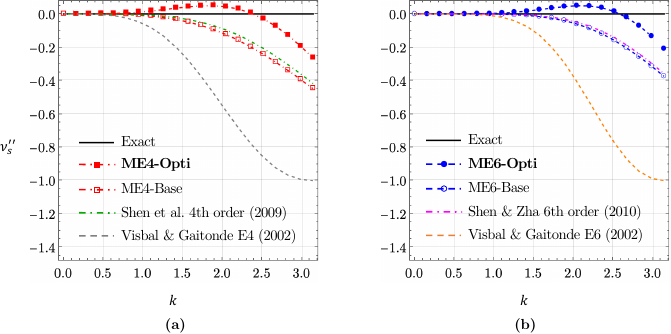}
    \caption{Spectral viscosity profiles of various 4th and 6th order straight derivative term discretization schemes.}
    \label{spectral_visc}
\end{figure}

From Fig. \ref{spectral_visc}, it is evident that as the wavenumber varies, so does the equivalent numerical dissipation. For example, when computing second derivatives in a solution field characterized by low-wavenumber features, the additional numerical dissipation is negligible, as the spectral viscosity at low wavenumbers (in $k \in [0,1]$) is nearly zero for all schemes shown in Fig. \ref{spectral_visc}. However, at higher wavenumbers, spectral viscosity can become negative, leading to anti-dissipation especially when $h$-ellipticity is absent. Notably, in the optimized schemes discussed, this anti-dissipation is minimized across the entire wavenumber spectrum $k \in [0,\pi]$, resulting in more accurate numerical modeling of the underlying physical phenomena. This has been demonstrated through various numerical experiments in sections 5.3, 5.4, and 5.5.

\subsection{Equivalent Reynolds number}
The equivalent Reynolds number ($Re_{\text{eq}}$), due to numerical dissipation in viscous discretization, is contingent on the spectral viscosity profile of the scheme. Given that spectral viscosity varies as a function of the local wavenumber, determined by the Fourier components in the discretized solution field, the equivalent Reynolds number also changes with wavenumber ($k$). By replacing $\nu$ with $\nu + \nu^{\prime \prime}_s$ in the Reynolds number formula, we obtain $Re_{\text{eq}} = \frac{L u}{\nu + \nu^{\prime \prime}_s}$. This adjustment leads to two distinct special scenarios as listed below:

\begin{equation}
    Re_{\text{eq}} = \frac{L u}{\nu + \nu^{\prime \prime}_s} = \frac{-L u k^2}{\nu k^{*}_{xx}} = \left\{
    \begin{array}{l}
    \text{Case-1: } k^{*}_{xx} = 0 \text{, arises due to odd-even mode at } k=\pi \text{,} \\
    \text{\phantom{Case-1: } \ \ \ \ \ \ \ \ \quad }Re_{\text{eq}} = \infty \text{ or } \frac{1}{Re_{\text{eq}}} = 0,\\
    \text{\phantom{Case-1: } \ \ \ \ \ \ \ \ \quad leads to pure inviscid phenomenon at } k=\pi. \\
    \text{Case-2: } k^{*}_{xx} = -k^2 \text{, exact modeling of the diffusion phenomenon.}
    \end{array}
    \right.
\end{equation}


\noindent \textbf{Case-1:} In this scenario, due to the odd-even mode at $k=\pi$, the spectral Reynolds number approaches infinity, effectively nullifying the viscous effects and leading to an inviscid model at the grid cutoff wavenumber. \\

\noindent \textbf{Case-2:} Contrary to case-1, this scenario leads to an accurate modeling of viscous effects. This is only possible in an ideal situation where the solution wavenumber being discretized is near zero or the grid is infinitely refined. Under these conditions, artificial numerical dissipation approaches to zero, allowing for an accurate representation of the diffusion effects. \\

As highlighted by Lamballais et al. \cite{lamballais2019implicit}, additional dissipation can be introduced at high wavenumbers to create a filtering effect, dissipating energy in these regions, similar to approaches used in several sub-grid scale large eddy simulation models \cite{lamballais2011straightforward, lamballais2021viscous}. The effectiveness of this strategy on the accuracy of shear-dominated tests, such as the doubly periodic shear layer and Kelvin-Helmholtz instabilities, remains unexplored and will be the subject of future research, as it is beyond the scope of the current discussion.

\begin{remark}
    Since all the schemes discussed in this work are central in nature, only the real part of the Fourier error analysis is non-zero (manifesting as a dissipation effect), while the imaginary part is always zero. However, in non-central second derivative schemes \cite{hoffman2018numerical} (not presented in this work), the error associated with the imaginary part will be non-zero, and it manifests as phase error and creates a spatial shift of Fourier modes corresponding to the discretized second derivative.
\end{remark}
\endgroup

\section{MP-OURS6 Inviscid flux discretization} \label{sec:OURS-6}
\begin{figure}[h!]
    \centering
    \includegraphics[width=\textwidth*7/8]{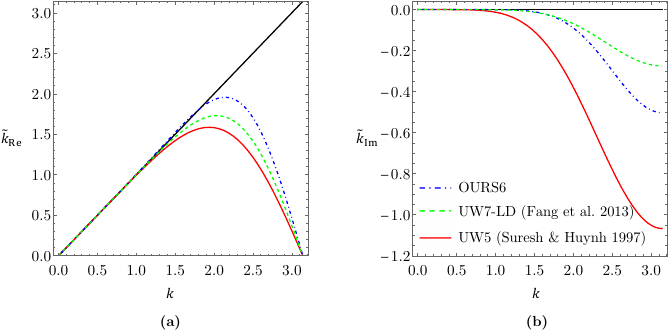}
    \caption{Spectral properties of OURS6 linear upwind interpolation scheme compared with the standard fifth order upwind \cite{suresh1997accurate,Shu1988} and seventh order bandwidth optimized low dissipation scheme of Fang et al. \cite{fang2013optimized}. (a) Real part (dispersion error) and (b) Imaginary part (dissipation error).}
    \label{fig:DsprDisp}
\end{figure}

The inviscid flux discretization employed in the supersonic jet case (Sec \ref{sec:SupersonicJet}) is briefly detailed here. Similar to the works of Li et al. \cite{li2013optimized}, Fang et al. \cite{fang2013optimized}, and Ahn \& Lee \cite{ahn2020modified} we modify the base reconstruction polynomial of the MP5 scheme to enhance the spectral resolution of the scheme. The original seven point stencil is being replaced with a eleven point stencil which was found to yield better results. A step-by-step procedure of the method is detailed below:

\begin{enumerate}
    \item Consider vector $\mathbf{P}=[\rho,u,v,w,p]$ to represent the array of primitive variables. Firstly, the left and right biased primitive variable states are estimated at the midpoint locations using a sixth order accurate eleven point stencil polynomial as follows:
    \begin{equation} \label{UL_poly1}
        \mathbf{P}_{i+\frac{1}{2}}^L \approx \sum_{p=-5}^4 c_{p} \mathbf{P}_{i+p}, \quad \mathbf{P}_{i+\frac{1}{2}}^R \approx \sum_{p=-4}^5 c_{-p} \mathbf{P}_{i+p}
    \end{equation}
    The coefficients used in the above expressions are: 
    
    \begin{equation*}
        c_{p}|_{p=-5 \text{ to } 4}=\left[\frac{3383}{3150000}, \frac{-3064}{196875}, \frac{32191}{393750}, \frac{-50632}{196875}, \frac{229667}{315000}, \frac{463193}{787500}, \frac{-128839}{787500}, \frac{33053}{787500}, \frac{-18073}{3150000}\right]
    \end{equation*}

    These coefficients have been obtained through an optimization process similar to the one presented in the current work (in Sec \ref{optimizedSchemes}). The spectral properties of the present scheme are shown in Fig. \ref{fig:DsprDisp}.
    

    \item Next, the interpolated states $\mathbf{P}_{i+\frac{1}{2}}^L$ and $\mathbf{P}_{i+\frac{1}{2}}^R$ are projected into characteristic space via a characteristic variable transformation. This involves multiplying the five interpolated primitive variables (obtained from previous step) with the corresponding left eigenvectors ($\overline{\mathbf{R}}_{\xi}^{-1}$) of the flux Jacobian matrix. The complete set of eigenvectors for curvilinear governing equations can be found in \cite{chandravamsi2023application}. The transformed values obtained after this process are denoted as $\mathbf{W}_{i+\frac{1}{2}}^L$ and $\mathbf{W}_{i+\frac{1}{2}}^R$.
    
    \begin{equation}
        \mathbf{W}_{i+\frac{1}{2}}^L=\overline{\mathbf{R}}_{\xi}^{-1} \mathbf{P}_{i+\frac{1}{2}}^L, \quad \mathbf{W}_{i+\frac{1}{2}}^R=\overline{\mathbf{R}}_{\xi}^{-1} \mathbf{P}_{i+\frac{1}{2}}^R
    \end{equation}
    
    \item To ensure local solution monotonicity, we apply the Monotonicity Preserving (MP) limiting operation of Suresh and Huynh \cite{suresh1997accurate} to the transformed variables. This step is crucial for preventing oscillations near sharp discontinuities and preserving the discontinuous features in the flow. While the initial unlimited upwind polynomial resolves the flow features in the smooth regions when the present limiting is not activated, the MP limiting operation ensures discontinuous feature preservation. Consequently, we obtain the reconstructed states $(\mathbf{W}_{i+\frac{1}{2}}^L)_{\text{limited}}$ and $(\mathbf{W}_{i+\frac{1}{2}}^R)_{\text{limited}}$ after this step.

    \item To convert $(\mathbf{W}_{i+\frac{1}{2}}^L)_{\text{limited}}$ and $(\mathbf{W}_{i+\frac{1}{2}}^R)_{\text{limited}}$ back into physical space, the reconstructed characteristic variables are multiplied with the corresponding right eigenvectors ($\overline{\mathbf{R}}_{\xi}$) of the flux Jacobian matrix \cite{chandravamsi2023application}.

    \begin{equation}
        (\mathbf{P}_{i+\frac{1}{2}}^L)_{\text{limited}}=\overline{\mathbf{R}}_{\xi}\left(\mathbf{W}_{i+1 / 2}^L\right)_{\text {limited}}, \quad (\mathbf{P}_{i+\frac{1}{2}}^R)_{\text{limited}}=\overline{\mathbf{R}}_{\xi}\left(\mathbf{W}_{i+1 / 2}^R\right)_{\text {limited}}
    \end{equation}

    \item The reconstructed primitive variable values $(\mathbf{P}_{i+\frac{1}{2}}^L)_{\text{limited}}$ and $(\mathbf{P}_{i+\frac{1}{2}}^R)_{\text{limited}}$ are used to compute the fluxes at each $i+\frac{1}{2}$ location using the HLLC (Harten-Lax-van Leer-Contact) approximate Riemann solver. Subsequently, the inviscid flux residual at each $i$ location is computed as follows:

    \begin{equation}
        \left(\frac{\partial \boldsymbol{\hat{F}}}{\partial \xi}\right)_i \approx \frac{\boldsymbol{\hat{F}}_{i+\frac{1}{2}}-\boldsymbol{\hat{F}}_{i-\frac{1}{2}}}{\Delta \xi}
    \end{equation}
\end{enumerate}

Named MP-OURS6, the method integrates the Monotonicity Preserving (MP) limiter for shock capturing and employs an Optimized sixth-order Upwind polynomial in the Reconstruction Scheme (OURS6).

%% file: _main.bbl
\begin{thebibliography}{10}
\expandafter\ifx\csname url\endcsname\relax
  \def\url#1{\texttt{#1}}\fi
\expandafter\ifx\csname urlprefix\endcsname\relax\def\urlprefix{URL }\fi
\expandafter\ifx\csname href\endcsname\relax
  \def\href#1#2{#2} \def\path#1{#1}\fi

\bibitem{wang2013high}
Z.~J. Wang, K.~Fidkowski, R.~Abgrall, F.~Bassi, D.~Caraeni, A.~Cary, H.~Deconinck, R.~Hartmann, K.~Hillewaert, H.~T. Huynh, et~al., High-order cfd methods: current status and perspective, International Journal for Numerical Methods in Fluids 72~(8) (2013) 811--845.

\bibitem{laizet2009high}
S.~Laizet, E.~Lamballais, High-order compact schemes for incompressible flows: A simple and efficient method with quasi-spectral accuracy, Journal of Computational Physics 228~(16) (2009) 5989--6015.

\bibitem{delorme2021application}
Y.~Delorme, R.~Stanly, S.~H. Frankel, D.~Greenblatt, Application of actuator line model for large eddy simulation of rotor noise control, Aerospace Science and Technology 108 (2021) 106405.

\bibitem{lele1992compact}
S.~K. Lele, Compact finite difference schemes with spectral-like resolution, Journal of Computational Physics 103~(1) (1992) 16--42.

\bibitem{bogey2004family}
C.~Bogey, C.~Bailly, A family of low dispersive and low dissipative explicit schemes for flow and noise computations, Journal of Computational physics 194~(1) (2004) 194--214.

\bibitem{liu2013new}
X.~Liu, S.~Zhang, H.~Zhang, C.-W. Shu, A new class of central compact schemes with spectral-like resolution i: Linear schemes, Journal of Computational Physics 248 (2013) 235--256.

\bibitem{shu1999high}
C.-W. Shu, High order eno and weno schemes for computational fluid dynamics, in: High-order methods for computational physics, Springer, 1999, pp. 439--582.

\bibitem{liu2009direct}
H.~Liu, J.~Yan, The direct discontinuous galerkin (ddg) methods for diffusion problems, SIAM Journal on Numerical Analysis 47~(1) (2009) 675--698.

\bibitem{huynh2009reconstruction}
H.~T. Huynh, A reconstruction approach to high-order schemnes including discontinuous galerkin for diffusion, in: 47th AIAA aerospace sciences meeting including the new horizons forum and aerospace exposition, 2009, p. 403.

\bibitem{gaitonde1998high}
D.~V. Gaitonde, M.~Visbal, High-order schemes for Navier-Stokes equations: algorithm and implementation into FDL3DI, Air Vehicles Directorte, Air Force Research Laboratory, Air Force Materiel~…, 1998.

\bibitem{shen2010large}
Y.~Shen, G.~Zha, Large eddy simulation using a new set of sixth order schemes for compressible viscous terms, Journal of Computational Physics 229~(22) (2010) 8296--8312.

\bibitem{hermeline2009finite}
F.~Hermeline, A finite volume method for approximating 3d diffusion operators on general meshes, Journal of computational Physics 228~(16) (2009) 5763--5786.

\bibitem{de2021high}
F.~De~Vanna, A.~Benato, F.~Picano, E.~Benini, High-order conservative formulation of viscous terms for variable viscosity flows, Acta Mechanica 232 (2021) 2115--2133.

\bibitem{tam1993dispersion}
C.~K. Tam, J.~C. Webb, Dispersion-relation-preserving finite difference schemes for computational acoustics, Journal of computational physics 107~(2) (1993) 262--281.

\bibitem{shen2009high}
Y.~Shen, G.~Zha, X.~Chen, High order conservative differencing for viscous terms and the application to vortex-induced vibration flows, Journal of Computational Physics 228~(22) (2009) 8283--8300.

\bibitem{wang2018accuracy}
N.~Wang, M.~Li, L.~Zhang, Accuracy analysis and improvement of viscous flux schemes in unstructured second-order finite-volume discretization, Chinese Journal of Theoretical and Applied Mechanics 50~(2018-3-527) (2018) 527.
\newblock \href {https://doi.org/10.6052/0459-1879-18-037} {\path{doi:10.6052/0459-1879-18-037}}.

\bibitem{nishikawa2011robust}
H.~Nishikawa, Robust and accurate viscous discretization via upwind scheme--i: Basic principle, Computers \& Fluids 49~(1) (2011) 62--86.

\bibitem{nishikawa2017effects}
H.~Nishikawa, Y.~Nakashima, N.~Watanabe, Effects of high-frequency damping on iterative convergence of implicit viscous solver, Journal of Computational Physics 348 (2017) 66--81.

\bibitem{nishikawa2010beyond}
H.~Nishikawa, Beyond interface gradient: a general principle for constructing diffusion schemes, in: 40th fluid dynamics conference and exhibit, 2010, p. 5093.

\bibitem{nishikawa2011two}
H.~Nishikawa, Two ways to extend diffusion schemes to navier-stokes schemes: Gradient formula or upwind flux, in: 20th AIAA Computational Fluid Dynamics Conference, 2011, p. 3044.

\bibitem{Visbal2002}
M.~R. Visbal, D.~V. Gaitonde, {On the Use of Higher-Order Finite-Difference Schemes on Curvilinear and Deforming Meshes}, Journal of Computational Physics 181~(1) (2002) 155--185.
\newblock \href {https://doi.org/10.1006/jcph.2002.7117} {\path{doi:10.1006/jcph.2002.7117}}.

\bibitem{fang2014investigation}
J.~Fang, Y.~Yao, Z.~Li, L.~Lu, Investigation of low-dissipation monotonicity-preserving scheme for direct numerical simulation of compressible turbulent flows, Computers \& Fluids 104 (2014) 55--72.

\bibitem{visbal1999high}
M.~R. Visbal, D.~V. Gaitonde, High-order-accurate methods for complex unsteady subsonic flows, AIAA journal 37~(10) (1999) 1231--1239.

\bibitem{sandham2002entropy}
N.~D. Sandham, Q.~Li, H.~C. Yee, Entropy splitting for high-order numerical simulation of compressible turbulence, Journal of Computational Physics 178~(2) (2002) 307--322.

\bibitem{Nishikawa2010}
H.~Nishikawa, {Beyond Interface Gradient: A General Principle for Constructing Diffusion Schemes}, 40th Fluid Dynamics Conference and Exhibit (2010).

\bibitem{chamarthi2022importance}
A.~S. Chamarthi, S.~Bokor, S.~H. Frankel, On the importance of high-frequency damping in high-order conservative finite-difference schemes for viscous fluxes, Journal of Computational Physics 460 (2022) 111195.

\bibitem{chamarthi2023role}
A.~S. Chamarthi, H.~Chandravamsi, N.~Hoffmann, S.~Bokor, S.~H. Frankel, On the role of spectral properties of viscous flux discretization for flow simulations on marginally resolved grids, Computers \& Fluids 251 (2023) 105742.

\bibitem{cheong2001grid}
C.~Cheong, S.~Lee, Grid-optimized dispersion-relation-preserving schemes on general geometries for computational aeroacoustics, Journal of Computational Physics 174~(1) (2001) 248--276.

\bibitem{fu2017targeted}
L.~Fu, X.~Y. Hu, N.~A. Adams, Targeted eno schemes with tailored resolution property for hyperbolic conservation laws, Journal of Computational Physics 349 (2017) 97--121.

\bibitem{lin2018optimization}
Y.~Lin, Y.~Chen, C.~Xu, X.~Deng, Optimization of a global seventh-order dissipative compact finite-difference scheme by a genetic algorithm, Applied Mathematics and Mechanics 39~(11) (2018) 1679--1690.

\bibitem{martin2006bandwidth}
M.~P. Mart{\'\i}n, E.~M. Taylor, M.~Wu, V.~G. Weirs, A bandwidth-optimized weno scheme for the effective direct numerical simulation of compressible turbulence, Journal of Computational Physics 220~(1) (2006) 270--289.

\bibitem{jin2018optimized}
Y.~Jin, F.~Liao, J.~Cai, Optimized low-dissipation and low-dispersion schemes for compressible flows, Journal of Computational Physics 371 (2018) 820--849.

\bibitem{li2013optimized}
X.-l. Li, Y.~Leng, Z.-w. He, Optimized sixth-order monotonicity-preserving scheme by nonlinear spectral analysis, International Journal for Numerical Methods in Fluids 73~(6) (2013) 560--577.

\bibitem{ashcroft2003optimized}
G.~Ashcroft, X.~Zhang, Optimized prefactored compact schemes, Journal of computational physics 190~(2) (2003) 459--477.

\bibitem{lamballais2011straightforward}
E.~Lamballais, V.~Fortun{\'e}, S.~Laizet, Straightforward high-order numerical dissipation via the viscous term for direct and large eddy simulation, Journal of Computational Physics 230~(9) (2011) 3270--3275.

\bibitem{lamballais2021viscous}
E.~Lamballais, R.~V. Cruz, R.~Perrin, Viscous and hyperviscous filtering for direct and large-eddy simulation, Journal of Computational Physics 431 (2021) 110115.

\bibitem{fu2021very}
L.~Fu, Very-high-order teno schemes with adaptive accuracy order and adaptive dissipation control, Computer Methods in Applied Mechanics and Engineering 387 (2021) 114193.

\bibitem{dairay2017numerical}
T.~Dairay, E.~Lamballais, S.~Laizet, J.~C. Vassilicos, Numerical dissipation vs. subgrid-scale modelling for large eddy simulation, Journal of Computational Physics 337 (2017) 252--274.

\bibitem{tsoutsanis2015comparison}
P.~Tsoutsanis, I.~W. Kokkinakis, L.~K{\"o}n{\"o}zsy, D.~Drikakis, R.~J. Williams, D.~L. Youngs, Comparison of structured-and unstructured-grid, compressible and incompressible methods using the vortex pairing problem, Computer Methods in Applied Mechanics and Engineering 293 (2015) 207--231.

\bibitem{drikakis2001spurious}
D.~Drikakis, P.~K. Smolarkiewicz, On spurious vortical structures, Journal of Computational physics 172~(1) (2001) 309--325.

\bibitem{mosedale2007assessment}
A.~Mosedale, D.~Drikakis, \href{https://doi.org/10.1115/1.2801374}{{Assessment of Very High Order of Accuracy in Implicit LES models}}, Journal of Fluids Engineering 129~(12) (2007) 1497--1503.
\newblock \href {http://arxiv.org/abs/https://asmedigitalcollection.asme.org/fluidsengineering/article-pdf/129/12/1497/5729131/1497\_1.pdf} {\path{arXiv:https://asmedigitalcollection.asme.org/fluidsengineering/article-pdf/129/12/1497/5729131/1497\_1.pdf}}, \href {https://doi.org/10.1115/1.2801374} {\path{doi:10.1115/1.2801374}}.
\newline\urlprefix\url{https://doi.org/10.1115/1.2801374}

\bibitem{thornber2008numerical}
B.~Thornber, D.~Drikakis, Numerical dissipation of upwind schemes in low mach flow, International journal for numerical methods in fluids 56~(8) (2008) 1535--1541.

\bibitem{thornber2008entropy}
B.~Thornber, D.~Drikakis, R.~J. Williams, D.~Youngs, On entropy generation and dissipation of kinetic energy in high-resolution shock-capturing schemes, Journal of Computational Physics 227~(10) (2008) 4853--4872.

\bibitem{minion1997performance}
M.~L. Minion, D.~L. Brown, Performance of under-resolved two-dimensional incompressible flow simulations, ii, Journal of Computational Physics 138~(2) (1997) 734--765.

\bibitem{Sutherland1893}
W.~Sutherland, \href{https://doi.org/10.1080/14786449308620508}{Lii. the viscosity of gases and molecular force}, The London, Edinburgh, and Dublin Philosophical Magazine and Journal of Science 36~(223) (1893) 507--531.
\newblock \href {http://arxiv.org/abs/https://doi.org/10.1080/14786449308620508} {\path{arXiv:https://doi.org/10.1080/14786449308620508}}, \href {https://doi.org/10.1080/14786449308620508} {\path{doi:10.1080/14786449308620508}}.
\newline\urlprefix\url{https://doi.org/10.1080/14786449308620508}

\bibitem{kundu2015fluid}
P.~K. Kundu, I.~M. Cohen, D.~R. Dowling, Fluid mechanics, Academic press, 2015.

\bibitem{zingg2000comparison}
D.~Zingg, S.~De~Rango, M.~Nemec, T.~Pulliam, Comparison of several spatial discretizations for the navier--stokes equations, Journal of computational Physics 160~(2) (2000) 683--704.

\bibitem{de1999aerodynamic}
S.~De~Rango, D.~Zingg, Aerodynamic computations using a higher-order algorithm, in: 37th Aerospace Sciences Meeting and Exhibit, 1999, p. 167.

\bibitem{chandravamsi2023application}
H.~Chandravamsi, A.~S. Chamarthi, N.~Hoffmann, S.~H. Frankel, On the application of gradient based reconstruction for flow simulations on generalized curvilinear and dynamic mesh domains, Computers \& Fluids 258 (2023) 105859.

\bibitem{gottlieb1998total}
S.~Gottlieb, C.-W. Shu, Total variation diminishing runge-kutta schemes, Mathematics of computation 67~(221) (1998) 73--85.

\bibitem{hindmarsh1984stability}
A.~Hindmarsh, P.~Gresho, D.~Griffiths, The stability of explicit euler time-integration for certain finite difference approximations of the multi-dimensional advection--diffusion equation, International journal for numerical methods in fluids 4~(9) (1984) 853--897.

\bibitem{bell1989second}
J.~B. Bell, P.~Colella, H.~M. Glaz, A second-order projection method for the incompressible navier-stokes equations, Journal of computational physics 85~(2) (1989) 257--283.

\bibitem{clausen2013entropically}
J.~R. Clausen, Entropically damped form of artificial compressibility for explicit simulation of incompressible flow, Physical Review E 87~(1) (2013) 013309.

\bibitem{achu2021entropically}
S.~Achu, N.~R. Vadlamani, Entropically damped artificial compressibility solver using higher order finite difference schemes on curvilinear and deforming meshes, in: AIAA Scitech 2021 Forum, 2021, p. 0634.

\bibitem{brown1995performance}
D.~L. Brown, Performance of under-resolved two-dimensional incompressible flow simulations, Journal of Computational Physics 122~(1) (1995) 165--183.

\bibitem{vadlamani2018distributed}
N.~R. Vadlamani, P.~G. Tucker, P.~Durbin, Distributed roughness effects on transitional and turbulent boundary layers, Flow, Turbulence and Combustion 100 (2018) 627--649.

\bibitem{delorme2017simple}
Y.~T. Delorme, K.~Puri, J.~Nordstrom, V.~Linders, S.~Dong, S.~H. Frankel, A simple and efficient incompressible navier--stokes solver for unsteady complex geometry flows on truncated domains, Computers \& Fluids 150 (2017) 84--94.

\bibitem{bogey2009shock}
C.~Bogey, N.~De~Cacqueray, C.~Bailly, A shock-capturing methodology based on adaptative spatial filtering for high-order non-linear computations, Journal of Computational Physics 228~(5) (2009) 1447--1465.

\bibitem{ryu2000magnetohydrodynamic}
D.~Ryu, T.~Jones, A.~Frank, The magnetohydrodynamic kelvin-helmholtz instability: A three-dimensional study of nonlinear evolution, The Astrophysical Journal 545~(1) (2000) 475.

\bibitem{san2015evaluation}
O.~San, K.~Kara, Evaluation of riemann flux solvers for weno reconstruction schemes: Kelvin--helmholtz instability, Computers \& Fluids 117 (2015) 24--41.

\bibitem{bogey2011finite}
C.~Bogey, N.~De~Cacqueray, C.~Bailly, Finite differences for coarse azimuthal discretization and for reduction of effective resolution near origin of cylindrical flow equations, Journal of Computational Physics 230~(4) (2011) 1134--1146.

\bibitem{quirk1997contribution}
J.~J. Quirk, A contribution to the great Riemann solver debate, Springer, 1997.

\bibitem{gallice2022entropy}
G.~Gallice, A.~Chan, R.~Loub{\`e}re, P.-H. Maire, Entropy stable and positivity preserving godunov-type schemes for multidimensional hyperbolic systems on unstructured grid, Journal of Computational Physics 468 (2022) 111493.

\bibitem{Jiang1995}
G.-S. Jiang, C.-W. Shu, {Efficient Implementation of Weighted ENO Schemes}, Journal of Computational Physics 126~(126) (1995) 202--228.

\bibitem{toro2009riemann}
E.~Toro, Riemann Solvers and Numerical Methods for Fluid Dynamics: A Practical Introduction, Springer Berlin Heidelberg, 2009.

\bibitem{kakumani2022use}
H.~C.~V. Kakumani, N.~R. Vadlamani, P.~G. Tucker, On the use of high order central difference schemes for differential equation based wall distance computations, Computers \& Fluids 248 (2022) 105666.

\bibitem{kumar2013weno}
G.~Kumar, S.~S. Girimaji, J.~Kerimo, Weno-enhanced gas-kinetic scheme for direct simulations of compressible transition and turbulence, Journal of Computational Physics 234 (2013) 499--523.

\bibitem{suresh1997accurate}
A.~Suresh, H.~Huynh, Accurate monotonicity-preserving schemes with runge-kutta time stepping, Journal of Computational Physics 136~(1) (1997) 83--99.

\bibitem{edgington2019aeroacoustic}
D.~Edgington-Mitchell, Aeroacoustic resonance and self-excitation in screeching and impinging supersonic jets--a review, International Journal of Aeroacoustics 18~(2-3) (2019) 118--188.

\bibitem{gojon2019antisymmetric}
R.~Gojon, E.~Gutmark, M.~Mihaescu, Antisymmetric oscillation modes in rectangular screeching jets, AIAA Journal 57~(8) (2019) 3422--3441.

\bibitem{zhang2002broadband}
Z.~Zhang, L.~Mongeau, S.~H. Frankel, Broadband sound generation by confined turbulent jets, The Journal of the Acoustical Society of America 112~(2) (2002) 677--689.

\bibitem{kakumani2023gpu}
H.~C.~V. Kakumani, A.~S. Chamarthi, N.~Hoffmann, S.~H. Frankel, Gpu-accelerated numerical study of temperature effects in choked under-expanded supersonic jets, in: AIAA SCITECH 2023 Forum, 2023, p. 0976.

\bibitem{zhao2000effects}
W.~Zhao, S.~H. Frankel, L.~Mongeau, Effects of spatial filtering on sound radiation from a subsonic axisymmetric jet, AIAA journal 38~(11) (2000) 2032--2039.

\bibitem{henderson2005experimental}
B.~Henderson, J.~Bridges, M.~Wernet, An experimental study of the oscillatory flow structure of tone-producing supersonic impinging jets, Journal of Fluid Mechanics 542 (2005) 115--137.

\bibitem{gojon2017flow}
R.~Gojon, C.~Bogey, Flow structure oscillations and tone production in underexpanded impinging round jets, AIAA Journal 55~(6) (2017) 1792--1805.

\bibitem{burns1998numerical}
J.~Burns, A.~Balogh, D.~Gilliam, V.~Shubov, Numerical stationary solutions for a viscous burgers' equation, Journal of Mathematical Systems Estimation and Control 8 (1998) 253--256.

\bibitem{mittal2011numerical}
R.~Mittal, G.~Arora, Numerical solution of the coupled viscous burgers’ equation, Communications in Nonlinear Science and Numerical Simulation 16~(3) (2011) 1304--1313.

\bibitem{tucker2003differential}
P.~Tucker, Differential equation-based wall distance computation for des and rans, Journal of Computational Physics 190~(1) (2003) 229--248.

\bibitem{lamballais2019implicit}
E.~Lamballais, T.~Dairay, S.~Laizet, J.~Vassilicos, Implicit/explicit spectral viscosity and large-scale sgs effects, in: Direct and Large-Eddy Simulation XI, Springer, 2019, pp. 107--113.

\bibitem{hoffman2018numerical}
J.~D. Hoffman, S.~Frankel, Numerical methods for engineers and scientists, CRC press, 2018.

\bibitem{fang2013optimized}
J.~Fang, Z.~Li, L.~Lu, An optimized low-dissipation monotonicity-preserving scheme for numerical simulations of high-speed turbulent flows, Journal of Scientific Computing 56~(1) (2013) 67--95.

\bibitem{ahn2020modified}
M.-H. Ahn, D.-J. Lee, Modified monotonicity preserving constraints for high-resolution optimized compact scheme, Journal of Scientific Computing 83~(2) (2020) 1--27.

\bibitem{Shu1988}
C.~W. Shu, S.~Osher, {Efficient implementation of essentially non-oscillatory shock-capturing schemes}, Journal of Computational Physics 77~(2) (1988) 439--471.
\newblock \href {https://doi.org/10.1016/0021-9991(88)90177-5} {\path{doi:10.1016/0021-9991(88)90177-5}}.

\end{thebibliography}
